\documentclass[10pt]{article}
\usepackage{fancyhdr}
\usepackage{extramarks}
\usepackage{amsmath}
\usepackage{amsthm}
\usepackage{amsfonts}
\usepackage{siunitx}
\usepackage{tikz}
\usepackage{algorithm}
\usepackage{algpseudocode}
\usepackage{multirow}
\usepackage[export]{adjustbox}
\usepackage{booktabs}
\usepackage{palatino}
\usepackage{graphicx}
\usepackage{subfigure}
\usepackage[colorlinks,linkcolor=black,anchorcolor=black,citecolor=black,urlcolor=blue]{hyperref}
\usepackage{amsmath,bm}
\usepackage{booktabs}
\usepackage{mathtools}
\usepackage{amssymb}
\usepackage{caption}
\usepackage{capt-of}
\usepackage{mciteplus}
\usepackage{cite}
\usepackage{mathrsfs}
\usepackage[title,titletoc,toc]{appendix}
\usepackage{xr}
\usepackage{parskip}
\usepackage{soul}
\usepackage{textcomp}
\usepackage[colaction]{multicol}
\usepackage[switch]{lineno}
\usepackage{lipsum}
\usepackage{etoolbox}
\usepackage{longtable}
\usepackage{array}
\usepackage{tablefootnote}
\usepackage{ragged2e}
\def\NoNumber#1{{\def\alglinenumber##1{}\State #1}\addtocounter{ALG@line}{-1}}
\newcommand{\algrule}[1][.2pt]{\par\vskip.5\baselineskip\hrule height #1\par\vskip.5\baselineskip}
\newcolumntype{C}[1]{>{\centering\arraybackslash}p{#1}}
\usetikzlibrary{automata,positioning}
\usetikzlibrary{automata,positioning}

\begin{document}

\title{PLPCA: Persistent Laplacian-enhanced PCA for Microarray Data Analysis}
\author{Sean Cottrell$^1$, Rui Wang$^1$, and Guo-Wei Wei$^{1,2,3}$\footnote{
		Corresponding author.		Email: weig@msu.edu} \\
$^1$ Department of Mathematics, \\
Michigan State University, MI 48824, USA.\\
$^2$ Department of Electrical and Computer Engineering,\\
Michigan State University, MI 48824, USA. \\
$^3$ Department of Biochemistry and Molecular Biology,\\
Michigan State University, MI 48824, USA.  
}
\date{\today} 

\maketitle

\begin{abstract}

Over the years, Principal Component Analysis (PCA) has served as the baseline approach for dimensionality reduction in gene expression data analysis. It primary objective is to identify a subset of disease-causing genes from a vast pool of thousands of genes. However, PCA possesses inherent limitations that hinder its interpretability, introduce classification ambiguity, and fail to capture complex geometric structures in the data. Although these limitations have been partially addressed in the literature by incorporating various regularizers such as graph Laplacian regularization, existing improved PCA methods still face challenges related to multiscale analysis and capturing higher-order interactions in the data. To address these challenges, we propose a novel approach called Persistent Laplacian-enhanced Principal Component Analysis (PLPCA). PLPCA amalgamates the advantages of earlier regularized PCA methods with persistent spectral graph theory, specifically persistent Laplacians derived from algebraic topology. In contrast to graph Laplacians, persistent Laplacians enable multiscale analysis through filtration and incorporate higher-order simplicial complexes to capture higher-order interactions in the data. We evaluate and validate the performance of PLPCA using benchmark microarray datasets that involve normal tissue samples and four different cancer tissues. Our extensive studies demonstrate that PLPCA outperforms all other state-of-the-art models for classification tasks after dimensionality reduction.


\end{abstract}
 \newpage

\setcounter{page}{1}
\renewcommand{\thepage}{{\arabic{page}}}


\clearpage

\section{Introduction}
Biological processes heavily depend on the different expression levels of genes over time. Thus, it is no surprise that analyzing gene expression data holds an important place in the field of biological and medical research, particularly in tasks such as identifying characteristic genes strongly correlated with various cancer types, as well as classifying tissue samples into cancerous and normal categories \cite{derisi1996use}.

In microarray analysis, messenger RNA (mRNA) molecules are collected from a tissue sample and converted into complementary DNA (cDNA). These cDNA molecules are subsequently labeled with a fluorescent dye and hybridized onto a microarray. The microarray is then scanned to measure the expression level of each gene. This process generates gene expression data, which represents the activity of each gene in a sample, typically at the RNA production level. The continuous advancements in microarray technology have led to the generation of large-scale gene expression datasets \cite{sarmah2011microarray}.

Gene expression data is commonly represented in matrix form, where rows correspond to genes and columns denote samples. Each matrix element indicates the expression level of a specific gene in a particular sample. Given the high dimensionality of gene expression data, often encompassing over 10,000 genes but merely a few hundred samples, considering all genes in a tumor classification analysis could introduce noise and notably augment computational complexity. Therefore, it is common practice to perform gene filtering or dimensionality reduction prior to applying classification methods in order to extract meaningful insights with reduced noise \cite{luecken2019current}.

There are several methods available in the literature for achieving effective dimensionality reduction, categorizable as either linear or non-linear based on the chosen distance metric. Notably, Principal Component Analysis (PCA) and Linear Discriminant Analysis (LDA) are examples of linear methods \cite{dunteman1989principal}. PCA remains widely employed, which considers the global Euclidean structure of the original data. On the other hand, LDA aims to find a linear combination of features that maximizes class separability and minimizes inter-class variance, particularly for multi-class classification problems \cite{xanthopoulos2013linear}. Non-linear methods divide into two subgroups: those preserving global pairwise distances and those maintaining local distances. Kernel PCA falls into the former category, while Laplacian Eigenmaps serve as an example of the latter (Laplacian Eigenmaps will be extensively discussed later) \cite{belkin2003laplacian, scholkopf1998nonlinear}. As the name suggests, Kernel PCA builds upon traditional PCA. While PCA may not perform well on datasets with complex algebraic/manifold structures that cannot be adequately represented in a linear space, Kernel PCA addresses this limitation by employing kernel functions in a reproducing kernel Hilbert space, thereby accommodating non-linearity.

While the methods mentioned have found broad applications in mathematical and statistical sciences, they also possess inherent limitations. In our specific context of gene expression data analysis, we aim to build upon recent advancements in PCA to mitigate the limitations intrinsic to advanced PCA techniques more effectively. By leveraging these improvements, we strive to enhance the analysis and interpretation of gene expression data. Although PCA is a widely-used procedure for dimensionality reduction, it has several associated weakness. These include a lack of interpretability of the principal components due to the dense loadings and issues with class ambiguities. Various methods have been proposed to address these issues. Most notably, Feng et al. proposed Supervised Discriminative Sparse PCA (SDSPCA) to include class information and sparse constraints by introducing a class label matrix and optimizing the $\text{L}_{2,1}$ norm \cite{feng2019supervised}. 

As we mentioned earlier, an additional desirable aspect of dimensionality reduction is the capability to identify low-dimensional structures that are embedded within higher-dimensional spaces. Graph theory has offered solutions to address this issue. In particular, Jiang et al. (2013) incorporated a graph Laplacian term into PCA (gLPCA) \cite{jiang2013graph}, while Zhang et al. (2022) combined this approach with SDSPCA to integrate structural information, interpretability, and class information \cite{zhang2022enhancing}. Additionally, they further developed a robust variant by employing the $\text{L}_{2,1}$ norm instead of the Frobenius norm in the loss function, and they proposed an iterative optimization algorithm \cite{jiang2013graph,zhang2022enhancing}. However, it is important to note that the graph Laplacian utilized in this method is only defined at a single scale and lacks topological considerations.

Eckmann introduced topological graphs using simplicial complexes, leading to the development of topological Laplacians on graphs \cite{eckmann1944harmonische}. Topological Laplacians generalize the traditional pairwise graph relations into many-body relations, and their kernel dimensions are identical to those of corresponding homological groups \cite{horak2013spectra}. This can be regarded  as a discrete generalization of the Hodge Laplacian on manifolds. Recently, we have introduced persistent Hodge Laplacians on manifolds \cite{chen2021evolutionary} and persistent combinatorial Laplacians on graphs \cite{wang2020persistent}. The latter is also known as persistent spectral graphs or persistent Laplacians (PLs) \cite{memoli2022persistent,wang2021hermes}. 
PLs can be viewed as a generalization of persistent homology \cite{zomorodian2004computing,edelsbrunner2008persistent,cang2017topologynet}. The fundamental idea of persistent homology is to represent data as a topological space, such as a simplicial complex. We can then use tools from algebraic topology to reveal the topological features of our data, such as holes and voids. Additionally, persistent homology employs filtration to perform a multiscale analysis of the data and thus creates a family of topological invariants to characterize data in a unique manner. Nevertheless, persistent homology cannot capture the homotopic shape evolution of data. PLs were designed to address this limitation \cite{wang2020persistent}. PLs have both harmonic spectra and non-harmonic spectra. The harmonic spectra with zero eigenvalues recover all the topological invariants from persistent homology, whereas the non-harmonic spectra reveal the homotopic shape evolution. PLs have been employed to facilitate machine learning-assisted protein engineering predictions \cite{qiu2023persistent}, accurately forecast future dominant SARS-CoV-2 variants BA.4/BA.5 \cite{chen2022persistent},  and predict protein-ligand binding affinity \cite{meng2021persistent}.

Our objective is to introduce PL-enhanced PCA theory (PLPCA). PLPCA can better capture multiscale geometrical structure information than standard graph regularization does. Specifically, PLs enhance our ability to recognize the stability of topological features in our data at multiple scales. This is achieved via filtration, which induces a sequence of simplicial complexes. We can study the spectra of its corresponding Laplacian matrix for each complex in the sequence to extract this topological and geometric information \cite{wang2021hermes}. Furthermore, we will validate our novel method and demonstrate its performance by microarray data analysis. 
 
Our work outlines the following: First, we delve into the mathematics behind PCA and its previous relevant improvements, such as sparseness, label information, and graph regularization. Next, we will discuss the tools from PL theory, which we believe may improve the effectiveness of dimensionality reduction. Later, we incorporate these tools to formulate two new PCA methods. The first one, denoted as pLPCA, is a simple persistent Laplacian-enabled PCA model. The persistent Laplacians introduce multiscale nonlinear geometric information to the gLPCA. The second method, denoted PLPCA, is persistent Laplacian-enhanced Robust Supervised Discriminative Sparse PCA. Lastly, we validate the proposed methods by a comparison of their results with those in the literature for tumor classifications after dimensionality reduction. Extensive results indicate that the proposed methods are the-state-of-art models for dimensionality reduction.

\section{Data Summary}
The benchmark datasets utilized in this analysis were obtained from The Cancer Genome Atlas, which is a project aimed at cataloging genetic mutations that contribute to cancer through genome sequencing \cite{cancer2012comprehensive, raphael2017integrated, cancer2015comprehensive, farshidfar2017integrative}. The Cancer Genome Atlas specifically focuses on 33 cancer types that fulfill the following criteria: poor prognosis, significant public health impact, and availability of samples.

In our study, we focus on two datasets: the Multisource dataset and the COAD dataset \cite{cancer2012comprehensive, raphael2017integrated, cancer2015comprehensive, farshidfar2017integrative}. The Multisource dataset comprises normal tissue samples and three different cancer types along with their corresponding gene expression data. The included cancer types are cholangiocarcinoma (CHOL), head and neck squamous cell carcinoma (HNSCC), and pancreatic adenocarcinoma (PAAD). Specifically, the CHOL dataset consists of 45 samples (9 normal tissue samples and 36 cancer samples), the HNSCC dataset contains 418 samples (20 normal and 398 cancer samples), and the PAAD dataset consists of 180 samples (4 normal and 176 cancer samples). Each dataset encompasses 20,502 genes.

Furthermore, the COAD dataset consists of 281 samples (19 normal samples and 262 colon adenocarcinoma samples) and also spans 20,502 genes. Once again, we emphasize the significant imbalance between the number of samples and the number of features in gene expression data, underscoring the importance of dimensionality reduction.

\begin{table}[H]
    \centering
    \caption{Datasets Summary.}
    \label{tab: comparison outlier}
    {
    \begin{tabular}{ c|cccc } 
    \toprule
    Dataset & Category & \# of Samples & \# of Features & Cancer Description  \\
    \midrule
    \multirow{4}{6em}{MultiSource} 
        & CHOL & 36 & 20502 & Cholangiocarcinoma \\ 
        & HNCC & 398 & 20502& Head and Neck Squamous Cell Carcinoma \\ 
        & PAAD & 176 & 20502 & Pancreatic Adenocarcinoma \\
        & Normal & 33 & 20502 & Normal Tissue \\
    \midrule
    \multirow{2}{6em}{COAD} 
        & COAD & 262 & 20502 & Colon Adenocarcinoma \\
        & Normal & 19 & 20502 & Normal Tissue \\
    \bottomrule 
    \end{tabular} 
    }
\end{table}

\section{Methods on Dimensionality Reduction}


\subsection{Principal Component Analysis} 

Recognizing the importance of dimensionality reduction for tumor classification given gene sequencing data, we formally introduce the notion of PCA. The purpose of PCA is to map $m$-dimensional data (with $n$ samples), $\mathbf{X} \in \mathbb{R}^{n \times m}$, into a $k$-dimensional space such that $k<m$. This is accomplished via computing the principal components, which can be used to perform a change of basis into a desired, lower dimensional space. The principal components are obtained via an eigendecomposition of the data covariance matrix, for which the principal components are eigenvectors. Equivalently, principal components are expressed as linear combinations of the original variables which explain the most variance. The goal is then to describe the maximal amount of variation in the original data in the first  fewer components \cite{jolliffe2016principal}. 

Mathematically, the optimal $k$-dimensional space is given by solving
\begin{equation}\label{eq:1}
    \underset{\textbf{U},\textbf{Q}}{\mathrm{min}} \lVert \textbf{X} - \textbf{UQ}^{T} \rVert^2_{\rm F},
\end{equation} 
where $\textbf{U} = \{\vec{u}_1, ..., \vec{u}_k\}$, $\mathbf{U} \in \mathbb{R}^{k \times m}$, represents the principal directions in order of explained variation and $\lVert {\bf X} \rVert^2_{\rm F}=\sqrt{\sum^n_{i=1}\sum^m_{j=1}|x_{ij}|^2}$ is  the Frobenius norm of ${\bf X}$. In classical PCA, we take $\textbf{U}$ to be orthogonal, $\textbf{U}^T\textbf{U} = \textbf{I}$, though we can also apply the orthogonality constraint to $\textbf{Q} = \{\vec{q}_1, ..., \vec{q}_n\}$, $\mathbf{Q} \in \mathbb{R}^{n \times k}$, which represents the projected data points in a new space. 

\subsection{Sparseness and Discriminative Information}

  PCA requires that the principal components be obtained via a linear combination of variables with non-zero weightings (called loadings.) In the context of gene selection, each variable would then represent a specific gene \cite{luecken2019current}. There is then an unnecessarily added layer of complexity by enforcing that the loadings be non-zero, as most of the genes would be irrelevant to our analysis and we may wish to focus on only a few. Thus, the interpretation of our principal components would be aided by the allowance of zero weights via the introduction of sparse PCA. The mathematical formulation of sparse PCA can take several forms, and the inclusion of an L$_{2,1}$ norm penalty term on the projected data matrix is the method chosen for solving SDSPCA \cite{feng2019supervised}. The L$_{2,1}$ norm is defined as $\| {\bf Q }\|_{2,1}=\sum^n_{i=1} \| {\bf q}_i\|_2 $. First, we calculate the L$_2$ norm of each row, and then compute the L$_1$ norm of row-based  L$_2$ norms. 

One can further build upon this by also introducing discriminative information to reduce class ambiguity. Supervised discriminative sparse PCA (SDSPCA) obtains principal components by introducing supervised label information as well as sparse constraints \cite{feng2019supervised}. This is realized via the following optimization formula:
\begin{equation}\label{eq:2}
\underset{\textbf{U,Q,A}}{\mathrm{min}} \lVert \textbf{X} - \textbf{UQ}^{T} \rVert^2_{F} + \alpha \lVert \textbf{Y} - \textbf{AQ}^{T} \rVert_{F}^{2} + \beta \lVert \textbf{Q} \rVert_{2,1}, \text{ s.t. } \textbf{Q}^T\textbf{Q} = \textbf{I}.
\end{equation} 
Here, $\alpha$ and $\beta$ are scale weights balancing the class label and sparse constraint terms. We arbitrarily initialize  matrix $\mathbf{A} \in \mathbb{R}^{c\times k}$, and $\textbf{Y} \in \mathbb{R}^{c\times n}$ represents the one-hot
coding class indicator matrix. Furthermore, $c$ represents the number of classes in the data. The class indicator matrix consists of 0's and 1's, with the position of element 1 in each column representing the class label. The matrix can be defined as follows, with $s_{i,j}$ representing class labels
\begin{equation}\label{eq:3}
\mathbf{Y}_{i,j} = 
    \begin{cases} 
      1, &  \text{if } s_{i,j} =i,\text{ } j = 1,...,n,\text{ } i = 1,...,c \\
      0, & \text{otherwise} 
   \end{cases}
\end{equation}
SDSPCA incorporates both label information and sparsity into PCA, with the second and third terms guaranteeing discriminative ability and interpretability, respectively. 

\subsection{Intrinsic Geometrical Structure}

While SDSPCA improves performance relative to traditional PCA, one still wish to capture and preserve the geometric structure of our gene sequence data during dimensionality reduction, motivating the introduction of graph regularization \cite{jiang2013graph}. 

Graph Laplacian-based embedding preserves local geometric relationships while maximizing the smoothness with respect to the intrinsic manifold of the data set in the low embedding space. Equivalently, one wish to construct a representation for the data sampled from a low-dimensional manifold embedded in a higher dimensional space. It has been shown that this can be accomplished with   graphs with pairwise edges, specifically the Laplacian operator \cite{belkin2001laplacian}. 

The core algorithm proceeds as follows. First, for $n$ points $\vec{x}_1,...,\vec{x}_n \in \mathbb{R}^m$ we construct a weighted graph $W$ with the set of nodes for each point, and the set of edges connecting neighboring points to one another. We put an edge between points if they are adjacent, which we can choose to determine according to K-Nearest-Neighbors (KNNs) or some distance threshold \cite{belkin2001laplacian}. While less geometrically intuitive, the KNN framework tends to be simpler \cite{cai2010graph}. We then weigh the edges via the Gaussian kernel and can obtain the following matrix:
\begin{equation}\label{eq:4}
\mathbf{W_{ij}} = 
    \begin{cases} 
    e^{-\lVert x_i - x_j\lVert^2/\eta}, \text{ if } x_i \text{ and } x_j \text{ are connected; } \\
    0,  \text{ otherwise}. \\
    \end{cases}
\end{equation} 
The matrix defined above, which we associate with our graph, is the adjacency matrix, and it encodes our connectivity information. Also note the introduction of a scale parameter $\eta \in \mathbb{R}$, which we must optimally choose. Alternatively, we can weight each edge $\mathbf{W_{ij}} = 1$ if points i and j are connected, but the choice of Gaussian Kernel weighting can be justified \cite{belkin2001laplacian}. Our goal can now be viewed as mapping the weighted connected graph to a line such that connected points stay close together. This means choosing $q_i \in \mathbb{R}$ to minimize:
\begin{equation}\label{eq:5}
    \sum_{i,j}\lVert q_i-q_j \lVert^2 \mathbf{W_{ij}}.
\end{equation}
It can be shown that this problem reduces to computing eigenvalues and eigenvectors for the generalized eigenvector problem:
\begin{equation}\label{eq:6}
    \mathbf{L}\vec{q} = \lambda \mathbf{D} \vec{q},
\end{equation} 
where $\mathbf{D}$ is a diagonal weight matrix with entries equaling the row sums of the adjacency matrix \cite{belkin2001laplacian}. $\mathbf{L} = \mathbf{D}-\mathbf{W}$ is the weighted Laplacian matrix. Let $\vec{q}_1,...,\vec{q}_n$ be the solutions of the eigenvector problem ordered according to their eigenvalues. The image of $\vec{x}_i$ under the embedding into the lower dimensional space $\mathbb{R}^k$ is then given by $\mathbf{Q} = (\vec{q}_1(i),...,\vec{q}_k(i)).$ Thus, we need to minimize:
\begin{equation}\label{eq:7}
    \sum_{i,j}\lVert \mathbf{Q}_i - \mathbf{Q}_j \lVert^2 \mathbf{W_{ij}} = \text{tr}(\mathbf{Q}^T\mathbf{LQ}).
\end{equation} 
Thus, given our data matrix $\mathbf{X}$ and weighted graph $W$, we seek a low dimensional representation that is regularized by the data manifold encoded in $W$. Because $\mathbf{Q}$ in PCA and Laplacian embedding serve the same purpose, we set them equal and combine Eqs. (\ref{eq:1}) and  (\ref{eq:7}), giving rise to graph Laplacian PCA (gLPCA), which is implemented according to the following optimization formula \cite{jiang2013graph}:
\begin{equation}\label{eq:8}
    \underset{\textbf{U,Q}}{\mathrm{min}} \lVert \textbf{X} - \textbf{UQ}^{T} \rVert^2_{F} + \gamma \text{Tr}(\textbf{Q}^{T} \textbf{L} \textbf{Q}), \text{  s.t. } \textbf{Q}^T\textbf{Q} = \textbf{I},
\end{equation}
where the $\gamma$ parameter scales the geometrical structure capture. Next, Zhang et al. combined this methodology with SDSPCA to incorporate sparseness, structural information, and discriminative information into one procedure. This new method, called Laplacian Supervised Discriminative Sparse PCA is solved for by combining Equations (7,2) \cite{zhang2022enhancing}:
\begin{equation}\label{eq:9}
    \underset{\textbf{U,Q,A}}{\mathrm{min}} \lVert \textbf{X} - \textbf{UQ}^{T} \rVert^2_{F} + \alpha \lVert \textbf{Y} - \textbf{AQ}^{T} \rVert_{F}^{2} + \beta \lVert \textbf{Q} \rVert_{2,1} + \gamma \text{Tr}(\textbf{Q}^{T} \textbf{L} \textbf{Q}), \text{  s.t. } \textbf{Q}^T\textbf{Q} = \textbf{I},
\end{equation} 

However, Zhang et al. also noted that the Frobenius norm regularization in LSDSPCA is  sensitive to outliers. For robustness, one can replace Frobenius norm regularization with $\text{L}_{2,1}$-norm regularization, which result in the robust Laplacian Supervised Discriminative Sparse PCA (RLSDSPCA) \cite{zhang2022enhancing}:
\begin{equation}\label{eq:10}
\underset{\textbf{U,Q,A}}{\mathrm{min}} \lVert \textbf{X} - \textbf{UQ}^{T} \rVert_{2,1} + \alpha \lVert \textbf{Y} - \textbf{AQ}^{T} \rVert_{F}^{2} + \beta \lVert \textbf{Q} \rVert_{2,1} + \gamma \text{Tr}(\textbf{Q}^{T} \textbf{L} \textbf{Q}), \text{ s.t. } \textbf{Q}^T\textbf{Q} = \textbf{I}.
\end{equation}
Having integrated robustness, interpretability, class information, and geometric structural information, we now turn to replacing the graph regularization with Persistent Spectral Graphs \cite{wang2020persistent} to introduce multiscale analysis. 

\subsection{Persistent Laplacians}
Motivated by the success of persistent homology and multiscale graphs in analyzing biomolecular data, we turn to persistent spectral graph theory to enhance our ability to capture the multiscale geometric structure \cite{nguyen2019agl}. Like persistent homology, persistent spectral graph theory tracks the birth and death of topological features of a dataset  as they change over scales \cite{memoli2022persistent, chen2021evolutionary}. We carry out this analysis by using the filtration procedure on our  dataset to construct a family of geometric configurations \cite{wang2020persistent}. We then can study the topological properties of each configuration by its corresponding Laplacian matrix. The topological persistence can be studied through multiple successive configurations.  

We first must introduce the notion of a simplex. A 0-simplex is a node, a 1-simplex is an edge, a 2-simplex is a triangle, a 3-simplex is a tetrahedron, and so on. Generally, we consider $q$-simplices which we label $\sigma_q$. A simplicial complex is a way of approximating a topological space by gluing together lower-dimensional simplices in a specific way. More formally, a simplicial complex $K$ is a collection of simplices such that: 
\begin{enumerate}\label{eq:11}
    \item[1)] If $\sigma_q \in \text{K}$ and $\sigma_p$ is a face of $\sigma_q$ then $\sigma_p \in $ K;
    \item[2)] The nonempty intersection of any two simplices is a face of both simplices. 
\end{enumerate}

A $q$-chain is then a formal sum of $q$-simplices in a simplicial complex $K$ with coefficients in $\mathbb{Z}_2$. The set of all $q$-chains has a basis of the set of $q$-simplicies in $K$. This set forms a finitely generated free Abelian group $C_q(K)$. We then define the boundary operator to be a group homomorphism that relates the chain groups, $\partial_q : C_q(K) \rightarrow C_{q-1}(K)$ \cite{wang2021hermes}. 

We denote the $q$-simplex by its vertices $v_i$: $\sigma_q = [ v_0, v_1, ..., v_q] $. The boundary operator is then defined as: 
\begin{equation}\label{eq:12}
    \partial_q \sigma_q := \sum_{i=0}^{q} (-1)^i \sigma_{q-1}^i,
\end{equation}
where $\sigma_{q-1}^i = [v_0,...,\hat{v_i},...,v_q]$ is the $(q-1)$ simplex with $v_i$ removed. The sequence of chain groups connected by boundary operators is then called a chain complex: 
\begin{center}
    $... \xrightarrow{\partial_{q+2}} C_{q+1}(K) \xrightarrow{\partial_{q+1}} C_{q}(K) \xrightarrow{\partial_q} ...$.
\end{center}

The chain complex associated with a simplicial complex defines the $q^{th}$ homology group $H_q = \text{ker} \partial_q / \text{Im} \partial_q$. The dimension of $H_q$ is then the $q^{th}$ Betti number $\beta_q$, which captures the number of $q$-dimensional holes in the simplicial complex. We can also define a dual chain complex through the adjoint operator of $\partial_q$ defined on the dual spaces $C^q(K) \cong C_q^*(K)$. The coboundary operator $\partial_q^*: C^{q-1}(K) \rightarrow C^q(K) $ is defined as: 
\begin{equation}
    \partial^* \omega^{q-1}(c_q) \equiv \omega^{q-1}(\partial c_q),
\end{equation}
where $\omega^{q-1} \in C^{q-1}(K)$ is a $(q-1)$ co-chain, or a  homomorphism mapping a chain to the coefficient group, and $c_q \in C_q(K)$ is a $q$-chain. The homology of the dual chain complex is referred to as the cohomology. Now we can define the q-combinatorial Laplacian operator, $\Delta_q : C^q(K) \rightarrow C^q(K)$ as: 
\begin{equation}
    \Delta_q := \partial_{q+1} \partial_{q+1}^* + \partial_q^* \partial_q.
\end{equation}

Now, denote the matrix representation of the $q$-boundary operator with respect to the standard basis for $C_q(K)$ and $C_{q-1}(K)$ as $\mathcal{B}_q$ and the matrix representation of the $q$-coboundary operator as $\mathcal{B}_q^T$. We then can define the matrix representation of the $q$th-order Laplacian operator as $\mathcal{L}_q$: 
\begin{equation}
    \mathcal{L}_q = \mathcal{B}_{q+1}\mathcal{B}_{q+1}^T + \mathcal{B}_q^T \mathcal{B}_q
\end{equation}

It is well known that $\beta_q$ is also the multiplicity of zero in the spectrum of the Laplacian matrix which corresponds to that simplicial complex (the harmonic spectrum). Specifically: 
\begin{align*}
    \beta_0 &= \text{Number of connected components in $K$} \\
    \beta_1 &= \text{Number of holes in $K$} \\
    \beta_2 &= \text{Number of two-dimensional voids in $K$}
\end{align*}
and so on. The non-harmonic spectrum also contains other topological and shape information. 

\begin{figure}[H]
	\centering
	\includegraphics[width = 0.7\textwidth]{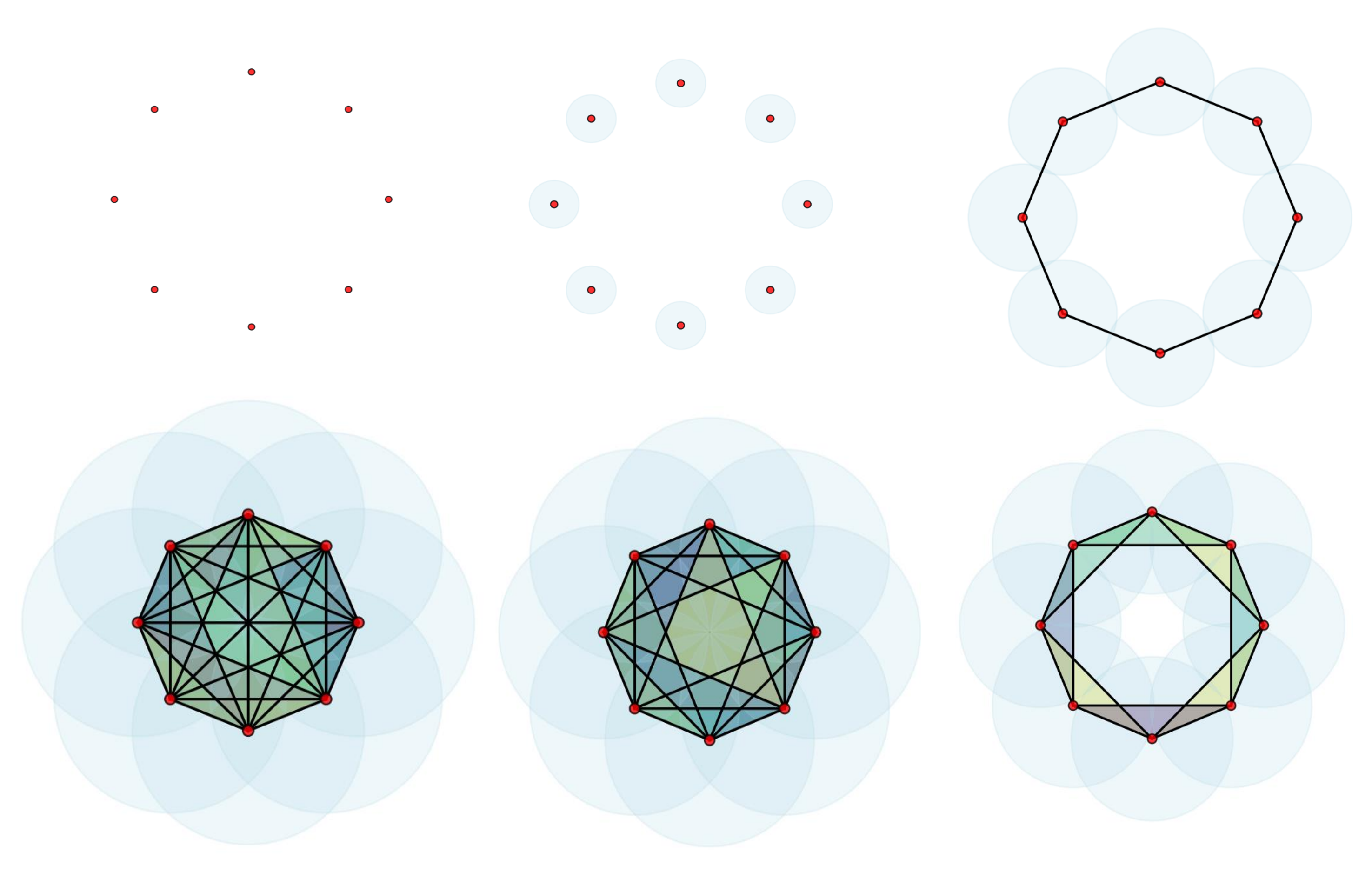}
	\caption{Illustration of the filtration of a point cloud by varying a distance threshold. The Vietoris-Rips complex is used.}
	\label{fig:filtration}
\end{figure}

However, a single simplicial complex offers very limited information to our understanding of the structure of our data. We then could consider the creation of a sequence of simplicial complexes which we induce by a filtration parameter: 
\begin{center}
    $\{ \emptyset \} = K_0 \subseteq K_1 \subseteq ... \subseteq K_p = K$
\end{center}
The filtration is illustrated in  Figure \ref{fig:filtration}. 

For each subcomplex $K_t$ we can denote its chain group to be $C_q(K_t)$, and the $q$-boundary operator $\partial_q^t: C_q(K_t) \rightarrow C_{q-1}(K_t)$. By convention, we define $C_q(K_t) = \{0\}$ for $q<0$ and the $q$-boundary operator to then be the zero map. We then have: 
\begin{equation}
    \partial_q^t \sigma_{q}= \sum_i^q(-1)^i\sigma_{q-1}^i, \forall \sigma_q \in K_t.
\end{equation}
Which is essentially the same construction as before. Likewise, the adjoint operator of $\partial_q^t$ is the coboundary operator $\partial_q^{t*}: C^{q-1}(K_t) \rightarrow C^q(K_t)$, which we regard as a map from $C_{q-1}(K_t) \rightarrow C_q(K_t)$ through the isomorphism between cochain and chain groups. We can then define a sequence of chain complexes.  

Next, we introduce persistence to the Laplacian spectra. Define the subset of $C_q^{t+p}$ whose boundary is in $C_{q-1}^t$ as $\mathbb{C}_q^{t,p}$,  assuming the natural inclusion map from $C_{q-1}^t \rightarrow C_{q-1}^{t+p}$.
\begin{equation}
    \mathbb{C}_q^{t,p} := \{ \beta \in C_q^{t+p} | \partial_q^{t,p} (\beta) \in C_{q-1}^t \}
\end{equation}
On this subset, one may define the $p$-persistent $q$-boundary operator denoted by $\hat{\partial}_q^{t,p}: \mathbb{C}_q^{t,p} \rightarrow C_{q-1}^t$ and corresponding adjoint operator $(\hat{\partial}_q^{t,p})^* : C_{q-1}^t \rightarrow  \mathbb{C}_q^{t,p}$, as before. The $q$-order $p$-persistent Laplacian operator is then $\Delta_q^{t,p}: C_q^t \rightarrow C_q^t$: 
\begin{equation}
    \Delta_q^{t,p} := \hat{\partial}_{q+1}^{t,p}(\hat{\partial}_{q+1}^{t,p})^* + (\hat{\partial}_q^{t})^* \hat{\partial}_q^{t}.
\end{equation}
And the matrix representation in simplicial basis is again: 
\begin{equation}
    \mathcal{L}_q^{t,p} := \mathcal{B}_{q+1}^{t,p}(\mathcal{B}_{q+1}^{t,p})^T + (\mathcal{B}_q^t)^T \mathcal{B}_q^t.
\end{equation}
We may again recognize the multiplicity of zero in the spectrum of $\mathcal{L}_q^{t,p}$ as the $q$th order $p$-persistent Betti number $\beta_q^{t,p}$ which counts the number of (independent) $q$-dimensional holes in $K_t$ that still exists in $K_{t+p}$. We can then see how the $q$th-order Laplacian is actually just a special case of the $q$th-order 0-persistent Laplacian at a simplicial complex $K_t$. In other words, the spectrum of $\mathcal{L}_q^{t,0}$ is simply associated with a snapshot of the filtration at some step $t$ \cite{wang2021hermes}. 

We can capture a more thorough view of the spatial features of our data by focusing on the 0-persistent Laplacian. Specifically, by inducing a family of subgraphs through the varying of a distance threshold $\epsilon$, as seen in Figure \ref{fig:filtration}. This is known as the Vietoris Rips Complex. The edges of the complex connect pairs of vertices that are within our distance threshold $\epsilon$, which we vary to construct our sequence of complexes. Alternatively, we connect vertices according to $K$-Nearest-Neighbors and then weight the edges according to some notion of distance, such as the Gaussian Kernel \cite{belkin2001laplacian}. We then filter out edges by increasing our $\epsilon$ value. In the next section we will see a convenient method for computation, and also a description of the PLPCA procedure. 

\subsection{PLPCA}

Now, the generation of the Vietoris Rips Complex can be achieved through implementing a filtration procedure on our weighted Laplacian matrix based on an increasing threshold. Observe:
\begin{equation}
\mathbf{L} = (l_{ij}), l_{ij} =
    \begin{cases}
    l_{ij}, i \neq j, i,j = 1,...,n \\
    l_{ii} = -\sum_{j=1}^{n}l_{ij}
    \end{cases} 
\end{equation}
For $i \neq j$, let $l_{\text{max}} = \text{max}(l_{ij})$, $l_{\text{min}} = \text{min}(l_{ij}), d = l_{\text{max}} - l_{\text{min}}$.
Set the $t^{th}$ Persistent Laplacian $\mathbf{L}^t, t = 1,...,p$: 
\begin{equation}
\mathbf{L}^t = (l_{ij}^t), l_{ij}^t = 
    \begin{cases}
    0, \text{ if } l_{ij} \leq (t/p)d + l_{\text{min}} \\
    -1, \text{otherwise}
    \end{cases}
\end{equation}
\begin{equation}
l_{ii}^t = -\sum_{j=1}^nl_{ij}^t
\end{equation}
This procedure results in the generation of a family of persistent Laplacians derived from our weighted graph Laplacian. However, due to the Gaussian Kernel weighting of the edges, it is more appropriate to filter values above the threshold rather than below. This is because more negative values indicate connections between data points that are closer together, which are the features we want to emphasize. To consolidate this family of graphs into a single term, we assign weights to each of them and then sum them together.

\begin{equation}\label{en:22}
\textbf{PL} := \sum_{t=1}^p\zeta_t\mathbf{L}^t.
\end{equation} 

Incorporating this new term into the gLPCA algorithm in place of the graph Laplacian should better retain geometrical structure information by emphasizing features that are persistent at multiple scales and providing a more thorough spatial view. The optimal space is now obtained via: 
\begin{equation}\label{en:23}
    \underset{\mathbf{U,Q}}{\mathrm{min}} \lVert \textbf{X} - \textbf{UQ}^{T} \rVert_{F}^2+ \gamma \text{Tr}(\textbf{Q}^{T} (\textbf{PL}) \textbf{Q}), \text{ s.t. } \textbf{Q}^T\textbf{Q} = \textbf{I}.
\end{equation}
We call this method pLPCA. We can then combine this procedure with RLSDSPCA to formulate the most optimal procedure, which, for convenience, we refer to simply as PLPCA, and which we solve for via:
\begin{equation}\label{en:24}
\underset{\mathbf{U,Q,A}}{\mathrm{min}} \lVert \textbf{X} - \textbf{UQ}^{T} \rVert_{2,1} + \alpha \lVert \textbf{Y} - \textbf{AQ}^{T} \rVert_{F}^{2} + \beta \lVert \textbf{Q} \rVert_{2,1} + \gamma \text{Tr}(\textbf{Q}^{T} (\textbf{PL}) \textbf{Q}), \text{ s.t. } \textbf{Q}^T\textbf{Q} = \textbf{I}.
\end{equation}

We should however recognize the issues with PLPCA regarding the inclusion of a class indicator matrix, most notably in the context of K-Means clustering via PCA \cite{ding2004k}. In this case, the label information would presumably not be known beforehand and therefore pLPCA would be the preferable method, although slightly less robust. 

\begin{figure}[H]
	\centering
    \includegraphics[ width = .5\textwidth]{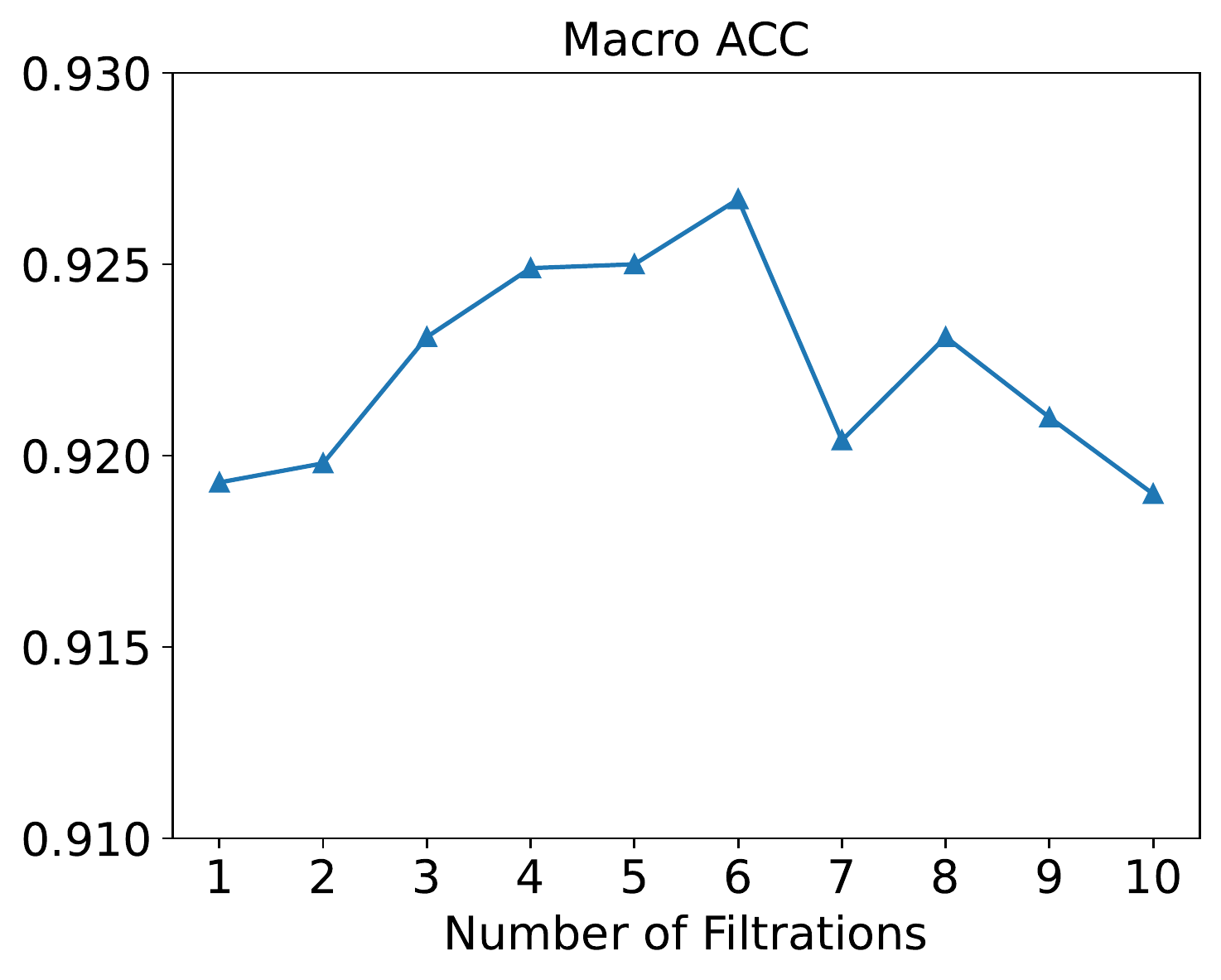}
	\caption{Effect of different numbers of filtrations on classification accuracy. The $x$-axis represents the number of flitrations and the $y$-axis represents the classification accuracy.  }
	\label{fig:p_value}
\end{figure}

\begin{figure}[H]
	\centering
  \includegraphics[width=0.7\linewidth]{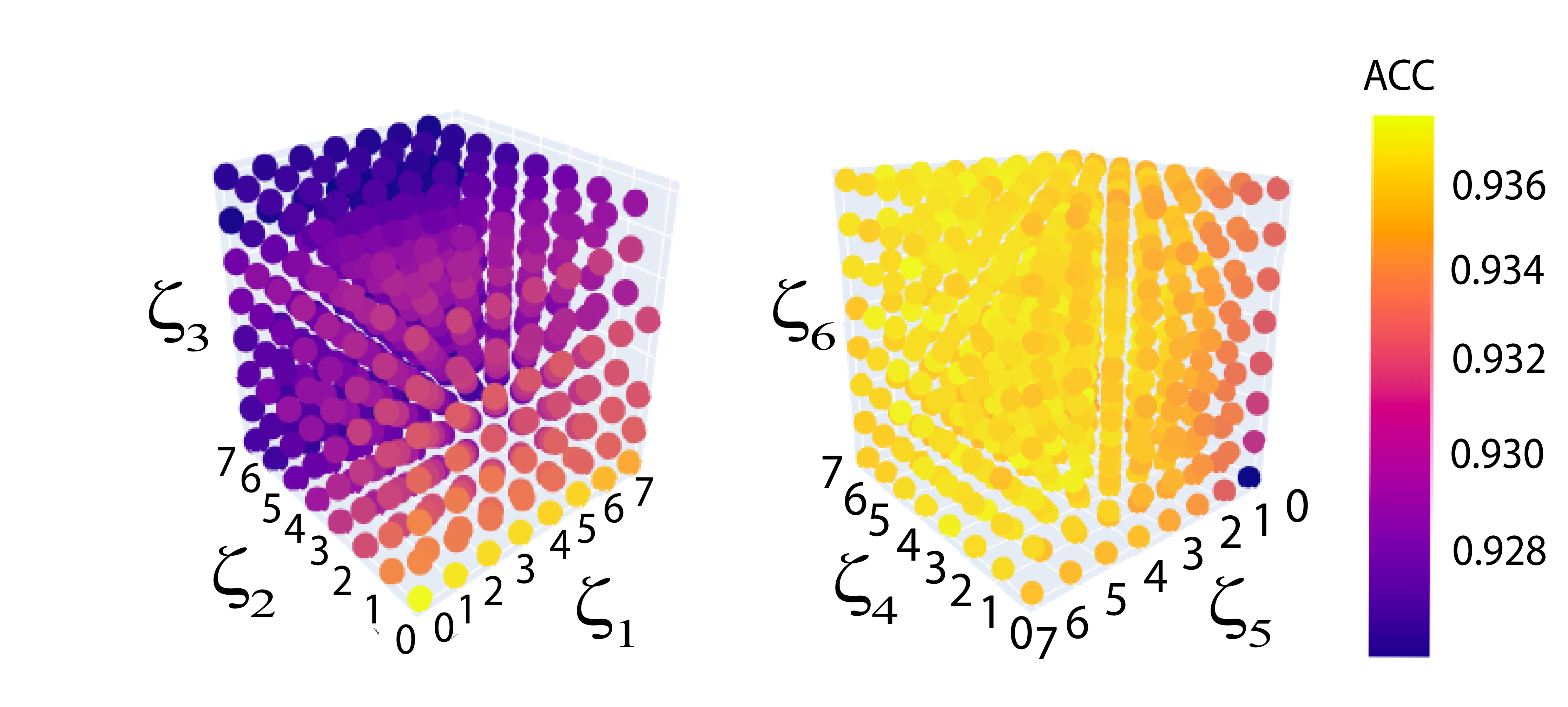}
  \caption{Effect of different close, middle, and long range connectivity weights on Mean Macro-ACC } \label{fig:ACC-z123}
 \end{figure}

Note also the introduction of several new hyper-parameters which we must optimize. Namely, each of the weights $\zeta_t$ as well as the number of sub-graphs must be chosen. We imposed the constraint $\sum_{t=1}^p \zeta_t = 1$ and performed grid search to understand whether we should favor long-range, middle-range, or close-range connectivity. Ultimately, we achieved the best results with six filtrations, or a six-scale scheme ($p=6$). Figure \ref{fig:p_value} shows the variations in Mean Macro-ACC as we vary the number of filtrations. For $p < 6$, there is not enough additional information incorporated to substantially improve performance, while for $p > 6$, too great a number of hyper-parameters to choose could have somewhat hurt performance. It is also possible that the additional filtrations hurt the performance when mapping into lower subspace dimensions, thereby hurting the performance on average.
For the MultiSource dataset, optimal results were achieved by emphasizing the long-range connections, while for the COAD dataset, the connectivities were all roughly equally weighted, with the close range and long range connectivity being slightly favored. Figure  \ref{fig:ACC-z123}  displays the optimization results for different combinations of weights in $\{ \zeta_t \}$. Higher values of $\zeta_t$ correspond to placing greater emphasis on the connectivity at that scale.
 
\begin{figure}[H]
	\centering
    \includegraphics[ width = .4\textwidth]{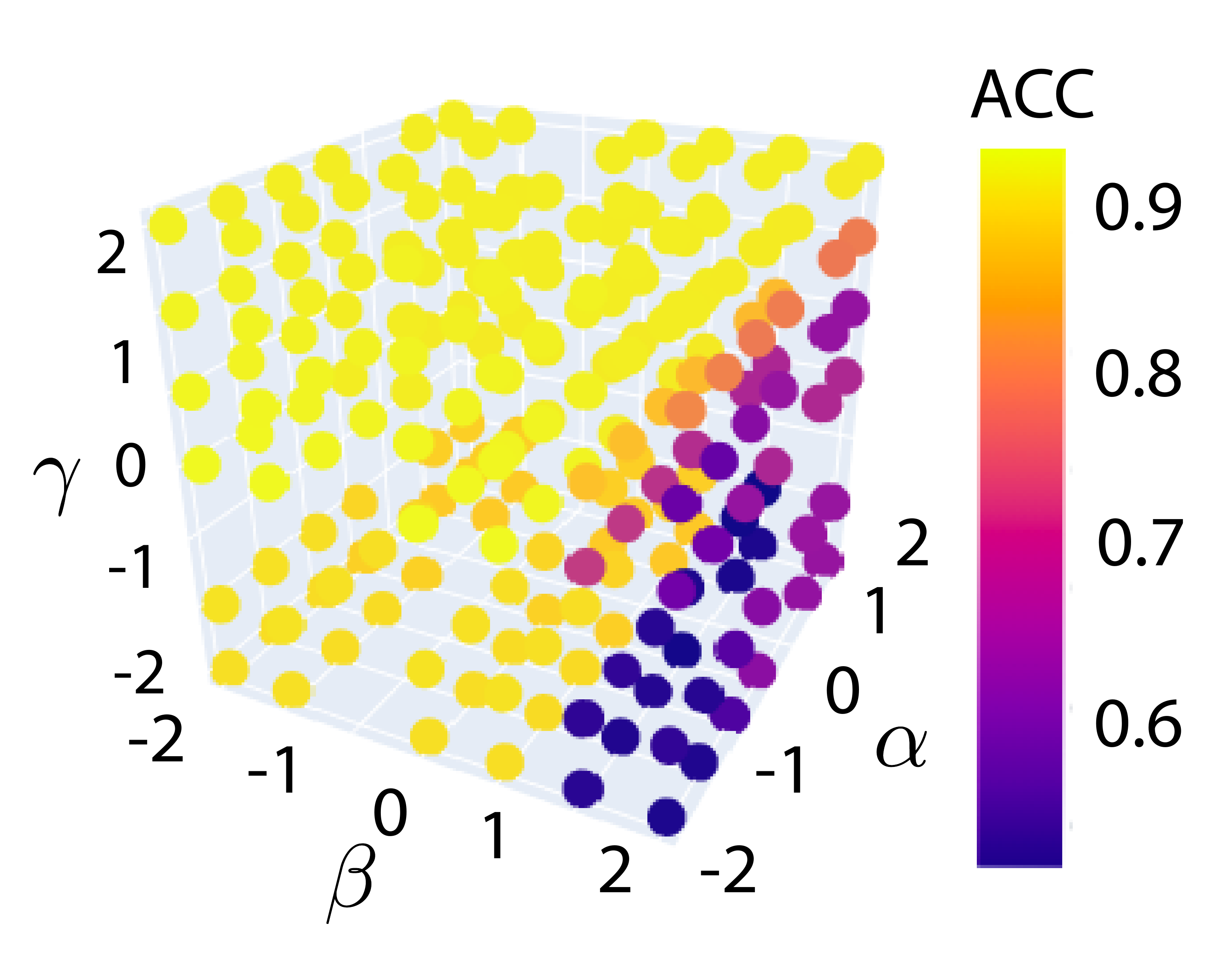}
	\caption{Effect of different weights on classification accuracy. The three coordinates represent each scale weight, color represents accuracy for each parameter combination. Axis ticks denote powers of 10.   }
	\label{fig:grid}
\end{figure}

We tested $\zeta_t$ values ranging from 0 to 10, and then scaled them to satisfy the constraint  $\sum_{t=1}^p \zeta_t = 1$. The difference was then factorized into the $\gamma$ parameter, which we test with each different combination of $\zeta_t$. For PLPCA on the COAD dataset, \{$\zeta_t\} = \{ 0.5,3,1,2,2,1 \}$. For pLPCA on the COAD dataset, \{$\zeta_t\} = \{2, 3, 0, 0, 2, 1\}$. For PLPCA on the MultiSource dataset, \{$\zeta_t\} = \{0.5, 0, 0, 3, 0, 6 \}$. For pLPCA on the MultiSource dataset, \{$\zeta_t\} = \{0.5, 0, 0, 3, 2, 6 \}$. After optimizing $\{\zeta_t \} \text{ and } \gamma$, we perform grid search again to revise our choices of $ \alpha \text{ and } \beta$ as shown in Figure \ref{fig:grid}.  
We revised our choice of parameters to $\alpha = 10^{-4}, \beta = 0.5, \text{and } \gamma = 10^{-1}$ for the MultiSource dataset and $\alpha = 10^{-5}, \beta = 0.5,\text{and } \gamma = 1000$ for the COAD dataset and pLPCA. For PLPCA, however, better results were obtained using $\gamma = 10^{-4}$ in both cases. We can now outline the optimization procedure behind our new method.

\begin{figure}[H]
	\centering
	\includegraphics[width = 0.8\textwidth]{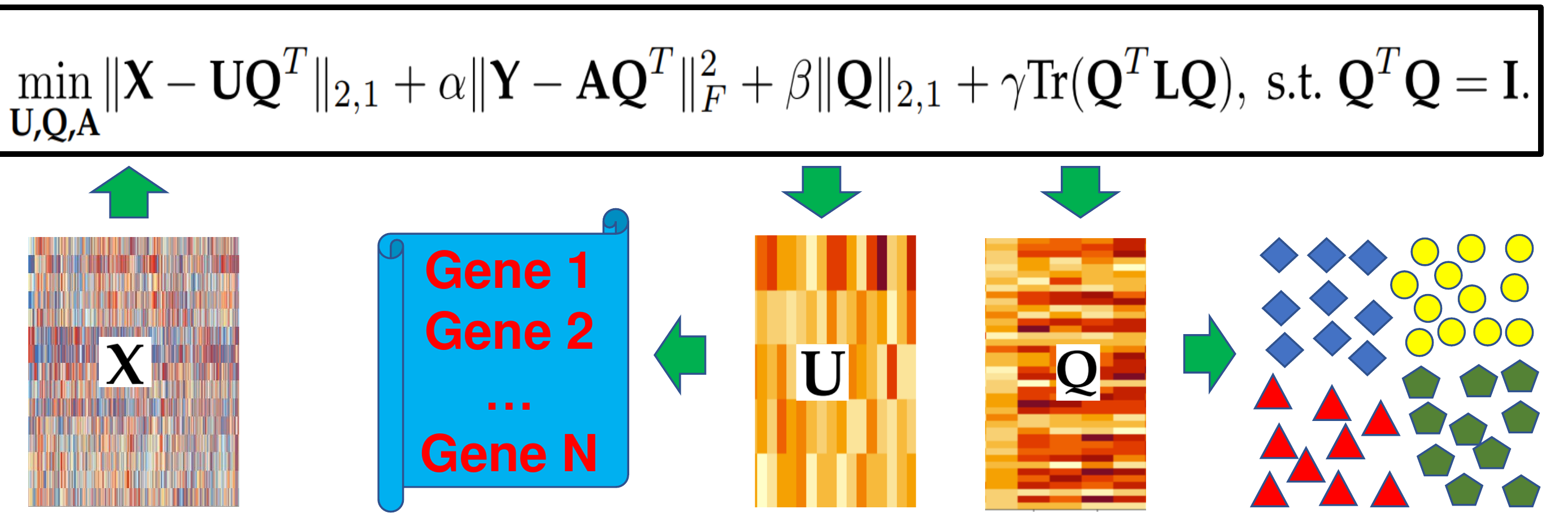}
	\caption{Outline of the PLPCA procedure for dimensionality reduction, feature selection, and classification.}
	\label{fig:procedure}
\end{figure}

Figure \ref{fig:procedure} provides an overview of the PLPCA framework for both tumor classification and characteristic gene selection \cite{zhang2022enhancing}. The process begins with the input gene expression matrix $\mathbf{X}$, and dimensionality reduction is performed using Eq. (\ref{en:24}). The resulting outputs are the projected data matrix $\mathbf{Q}$ and the principal directions matrix $\mathbf{U}$. The projected data matrix can then be utilized for tumor classification. Previous studies have also utilized the projected data matrix for feature selection, as it contains valuable information about each gene's contribution to the overall variance of the data \cite{zhang2022enhancing}.

The optimization of PLPCA can be performed using the alternating direction method of multipliers (ADMM) algorithm. ADMM is a variant of the augmented Lagrangian method, which is employed to solve constrained optimization problems \cite{zhang2022enhancing}. The augmented Lagrangian method transforms constrained optimization problems into a series of unconstrained problems by introducing a penalty term to the cost function. It also incorporates an additional term resembling a Lagrangian multiplier \cite{fortin2000augmented}. The penalty function approach solves this problem iteratively by updating each parameter at each step \cite{echebest2016convergence}.

ADMM, meanwhile, is a method which uses partial updates for dual variables \cite{boyd2011distributed}. Consider the generic problem:
\begin{equation}
    \underset{\mathbf{Q}}{\mathrm{min}}f(\mathbf{Q}) + g(\mathbf{U})
\end{equation}
Which is equivalent to: 
\begin{equation}
    \underset{\mathbf{Q,U}}{\mathrm{min}}f(\mathbf{Q}) + g(\mathbf{U}) \text{ subject to } \mathbf{QQ}^T = \mathbf{I} 
\end{equation}
The ADMM technique allows us to approximately solve this problem by first solving for $\mathbf{Q}$ with $\mathbf{U}$ fixed via the augmented Lagrangian method, and then vice versa. Specific to our problem, this approach can be taken for approximately solving for the optimal $\mathbf{Q, U, A, E, G,} \text{ and } \mathbf{C}$ matrices in our algorithm. This is the method chosen for optimizing RLSDSPCA, and we implemented it as well for PLPCA \cite{zhang2022enhancing}. 

The optimization procedure for RLSDSPCA is described in detail in the paper by Zhang et al. \cite{zhang2022enhancing}. Their study demonstrated that the RLSDSPCA objective function exhibits a monotonically decreasing trend with each iteration. The same proof applies to PLPCA as well, with the substitution of the Persistent Laplacian for the Laplacian term. In Algorithm \ref{Algo}, we present a summary of the updated optimization algorithm for PLPCA.

After obtaining the optimal dimensionality reduction, the classification of cancerous tumors begins by normalizing the gene expression data and randomly splitting it into training and testing sets. The testing set accounts for 20\% of the data, while the training set constitutes the remaining 80\%. To mitigate the impact of data distribution, we employed a 5-fold cross-validation approach. The classification accuracy was calculated as the average performance over five repetitions. For classification purposes, we utilized the K-Nearest-Neighbors algorithm \cite{peterson2009k}. The mean accuracy of the classification was recorded for subspace dimensions ranging from \{100, 95, ..., 5, 1\}.  

\begin{algorithm}[ht]
    \caption{PLPCA procedure}\label{Algo}
    \hspace*{\algorithmicindent} \textbf{Input:} 
    \begin{algorithmic}
        \State Data matrix $\textbf{X} \in \mathbb{R}^{n\times m}$; OHE Label Matrix $\textbf{Y} \in \mathbb{R}^{c\times n}$;
        \State Weight Parameters: $\alpha, \beta, \gamma$; Number of Subspace Dimensions: k;
        \State Convergence Parameter: $\Theta$, Number of Iterations: \textit{MaxIter};
        \State Weight Parameters $\zeta_i$ ($i = 1,...,p$); Number of Subgraphs: $p$
    \algrule
    \end{algorithmic}
    \hspace*{\algorithmicindent} \textbf{Output:}
    \begin{algorithmic}
        \State Principal Directions Matrix \textbf{U} and Projected Data Matrix \textbf{Q}
    \end{algorithmic}
    \algrule 
    \hspace*{\algorithmicindent} \textbf{Initialize:}
    \begin{algorithmic}[1]
        \State Initialize Matrices \textbf{G, E, C} to identity matrix;
        \State Randomly intialize matrices \textbf{A}, $\textbf{Q}_1$ (An auxillary matrix to check convergence);
        \State Construct Weighted Adjacency Matrix \textbf{W} according to K-Nearest-Neighbors;
        \State Compute Weighted Laplacian \textbf{L};
        \State Compute family of subgraphs $\textbf{L}^t$;
        \State Construct Persistent Laplacian \textbf{PL};
        \State Initialize $\mu = 1$ 
        \algrule 
        
        \NoNumber{\textbf{ADMM}}
        \For{i=1 to \textit{MaxIter}}
            \State Compute \textbf{Q}
            \State Compute \textbf{U}
            \State Compute \textbf{A}
            \State Compute \textbf{E}
            \State Compute \textbf{G}
            \State Compute \textbf{C}
            \State Compute $\mu$
            \State Check the Convergence Condition:
            \If{$i>1$ and $\lVert \textbf{Q} - \textbf{Q}_1 \rVert _{2,1} < \Theta$}
                \State Break
            \EndIf 
        \EndFor
        \State Let $\textbf{Q}_1 = \textbf{Q}$
    \end{algorithmic}
\end{algorithm}

\section{Applications in  Tumor Classification Problems: Results and Discussion}

\subsection{Evaluation Metrics}

 Our study demonstrates that the incorporation of persistent graph regularization enhances classification performance after dimensionality reduction, surpassing the achievements of other state-of-the-art methods such as RLSDSPCA. We summarize the outcomes of our analysis conducted on the COAD and MultiSource datasets, sourced from the Cancer Genome Atlas, as well as several simulated Outlier datasets obtained from the RLSDSPCA GitHub repository \cite{zhang2022enhancing}.

Next, we discuss the evaluation metrics used to measure performance. While accuracy is a commonly employed metric, in the context of cancer diagnosis, additional emphasis is often placed on the F1-score. The F1-score represents the harmonic mean of Precision and Recall, providing a balanced assessment of a classifier's performance.

\begin{align}
    \text{Recall} &= \text{True Positive} / \text{(True Positive + False Negative)} \\
    \text{Precision} &= \text{True Positive} / \text{(True Positive + False Positive)} \\
    \text{F1-Score} &= 2 \times (\text{Precision} \times \text{Recall})/(\text{Precision}+\text{Recall})
\end{align}

By considering the cost of False Negatives in our classification task, we acknowledge the significance of accurately identifying cases of a potentially life-threatening disease. The F1-score is particularly relevant in this context as it takes into account both precision and recall, making it more robust to class imbalances within the data. In our case, there are noticeable imbalances, particularly in the MultiSource dataset.

Given that our data consists of multiple categories, the evaluation criterion we employ is the mean of each category indicator. This evaluation approach is commonly known as a Macro metric, where performance measures are calculated for each category individually and then averaged to obtain an overall score.


\begin{align}
    \text{Macro-Recall} &= 1/c \times \sum_{i=1}^{c} \text{Recall}_i    \\
    \text{Macro-Precision} &= 1/c \times \sum_{i=1}^{c} \text{Precision}_i   \\
    \text{Macro-F1} &= 2 \times (\text{Macro-Precision} \times \text{Macro-Recall})/(\text{Macro-Precision}+\text{Macro-Recall})
\end{align}

To enhance visualization, Residue Similarity (R-S) scores can be computed \cite{hozumi2022ccp}. Traditional visualization techniques often involve reducing the data to two or three dimensions, which may result in the loss of structure and integrity in multiclass data. R-S plots were introduced as a method to visualize results while better preserving the underlying structure of the data.

An R-S plot consists of two main components: the residue score and the similarity score. The residue score is calculated as the sum of distances between classes, capturing the dissimilarity between them. On the other hand, the similarity score represents the average similarity within each class, indicating the degree of similarity between instances belonging to the same class. By considering both scores, R-S plots provide a comprehensive representation of the data's structure in a visualization.

Given data of the form $ \{(\vec{x}_m, y_m) | \vec{x}_m \in \mathbb{R}^N, y_m \in \mathbb{Z}_l \}_{m=1}^M$, we have $y_m$ representing the class label of our $m$th data point $\vec{x}_m \in \mathbf{X}$. Say that our data has $N$ samples, $M$ features, and $L$ classes. We can then partition our dataset $\mathbf{X}$ into subsets containing each of the classes by taking $\mathcal{C}_l = \{ \vec{x}_m \in \mathbf{X} | y_m = l\}$. For each class $l$ we then define the residue score as follows: 

\begin{equation}
    R_m := R(\vec{x}_m) = \frac{1}{R_{\rm max}} \sum_{\vec{x}_j \notin \mathcal{C}_l} \lVert \vec{x}_m - \vec{x}_j
 \rVert,
\end{equation} 

where $\lVert \cdot \rVert$ denotes the Euclidean distance between vectors and $R_{\rm max}$ is the maximal residue score for that subset. The similarity score, meanwhile, is given as: 

\begin{equation}
    S_m := S(\vec{x}_m) = \frac{1}{\lvert \mathcal{C}_l \rvert} \sum_{\vec{x}_j \in \mathcal{C}_l}(1 - \frac{\lVert \vec{x}_m - \vec{x}_j \rVert} {d_{\rm max}}), 
\end{equation} 

where $d_{\rm max}$ is the maximal pairwise distance of the dataset. For constructing R-S plots, we then take $R(\vec{x})$ to be the ${x}$-axis and $S(\vec{x})$ to be the $ {y}$-axis. 

\subsection{Comparison of pLPCA and gLPCA}
The classification of benchmark tumor datasets provides an opportunity to evaluate the performance of our method in comparison to other state-of-the-art approaches. Firstly, we validate our claim that persistent Laplacian-based regularization surpasses graph Laplacian regularization by comparing pLPCA and gLPCA. Following this validation, we proceed to compare our PLPCA with RLSDSPCA, which has demonstrated the best performance in the existing literature \cite{zhang2022enhancing}.

To summarize the comparison between pLPCA and gLPCA on the COAD dataset, please refer to Table \ref{tab: comparison RL vs RPL}.  Note that the gLPCA results reported in by  the early work \cite{zhang2022enhancing} do  not appear to be reasonable because they are higher than those of the improved model 
gLSPCA, which should not happen according to Feng et al. \cite{feng2019pca}  Therefore, we have reproduced these results for gLPCA using the parameters specified by Zhang el at. \cite{zhang2022enhancing}  Our results are listed in Table \ref{tab: comparison RL vs RPL} for a comparison. 

\begin{table}[H] 
    \centering
    \caption{Comparison of pLPCA and gLPCA performance on the COAD dataset.}
    \label{tab: comparison RL vs RPL}
    {
    \begin{tabular}{ cccccc } 
    \toprule
    Method & Mean ACC & Mean Macro-REC & Mean Macro-PRE & Mean Macro-F1 & Macro-AUC  \\
    \midrule 
		 \multirow{1}{6em}{gLPCA  \cite{jiang2013graph,zhang2022enhancing}}&0.9777 & 0.9463 & 0.9138 & 0.9249 & 0.9470\\  
    \multirow{1}{6em}{gLPCA* \cite{jiang2013graph}} &0.9756   &0.9429  &0.8841  &0.9002  &0.9429 \\ 
    \multirow{1}{6em}{pLPCA} & \textbf{0.9788} & \textbf{0.9450} & \textbf{0.8996} & \textbf{0.9115} & \textbf{0.9450}  \\ 
    \bottomrule
    \end{tabular} 
		*Reproduced in the present work
    }
\end{table}

This table clearly demonstrates that pLPCA outperforms gLPCA in all performance metrics, highlighting the superior ability of persistent spectral graphs to retain topological and geometrical information during dimensionality reduction. Specifically, we observe an improvement in Mean Accuracy from 0.9756 to 0.9788. Similarly, the Macro-F1 score has improved from 0.9002 to 0.9115. These results indicate that pLPCA not only achieves higher overall classification accuracy but also reduces the number of False Negatives.

To further validate the performance of pLPCA, we conducted tests on the MultiSource dataset, and the results are presented in Table \ref{tab: comparison PL vs gL}. 

\begin{table}[H]
    \centering
    \caption{Comparison of pLPCA and gLPCA performance on the MultiSource dataset.}
    \label{tab: comparison PL vs gL}
    {
    \begin{tabular}{ cccccc } 
    \toprule
    Method & Mean ACC & Mean Macro-REC & Mean Macro-PRE & Mean Macro-F1 & Macro-AUC  \\
    \midrule
    \multirow{1}{6em}{gLPCA \cite{jiang2013graph,zhang2022enhancing}} & 0.9139 & 0.8147 & 0.8768  & 0.8316 &  0.8909\\
    \multirow{1}{6em}{pLPCA} & \textbf{0.9267} & \textbf{0.8318} & \textbf{0.8857} & \textbf{0.8471} & \textbf{0.8991} \\
    \bottomrule
    \end{tabular} 		   
    }
\end{table}

Once again, it is important to highlight the significant improvement in performance across all five evaluation metrics achieved by incorporating persistent Laplacian-based regularization instead of graph Laplacian. The mean accuracy has shown a substantial improvement of 1.28\%, while the Macro-F1 score has improved by 1.55\%.

\begin{figure}[H]
\minipage{0.5\textwidth}
  \includegraphics[width=\linewidth]{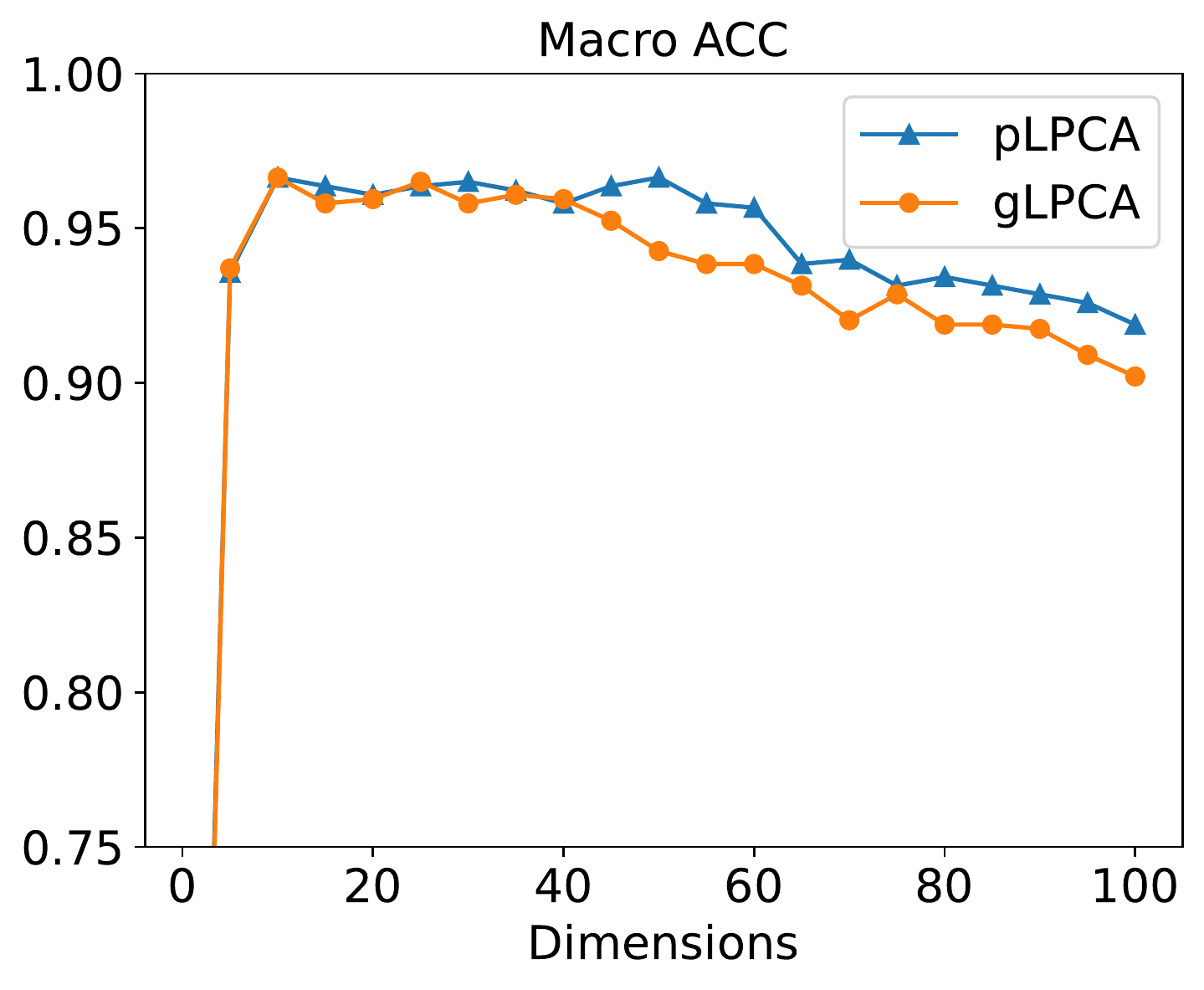}
  \caption{Macro ACC score on the MultiSource dataset across different reduced dimensions (gLPCA vs. pLPCA).}\label{fig:ACC-pLPCA}
\endminipage
\hspace{.5cm}
\minipage{0.5\textwidth}
  \includegraphics[width=\linewidth]{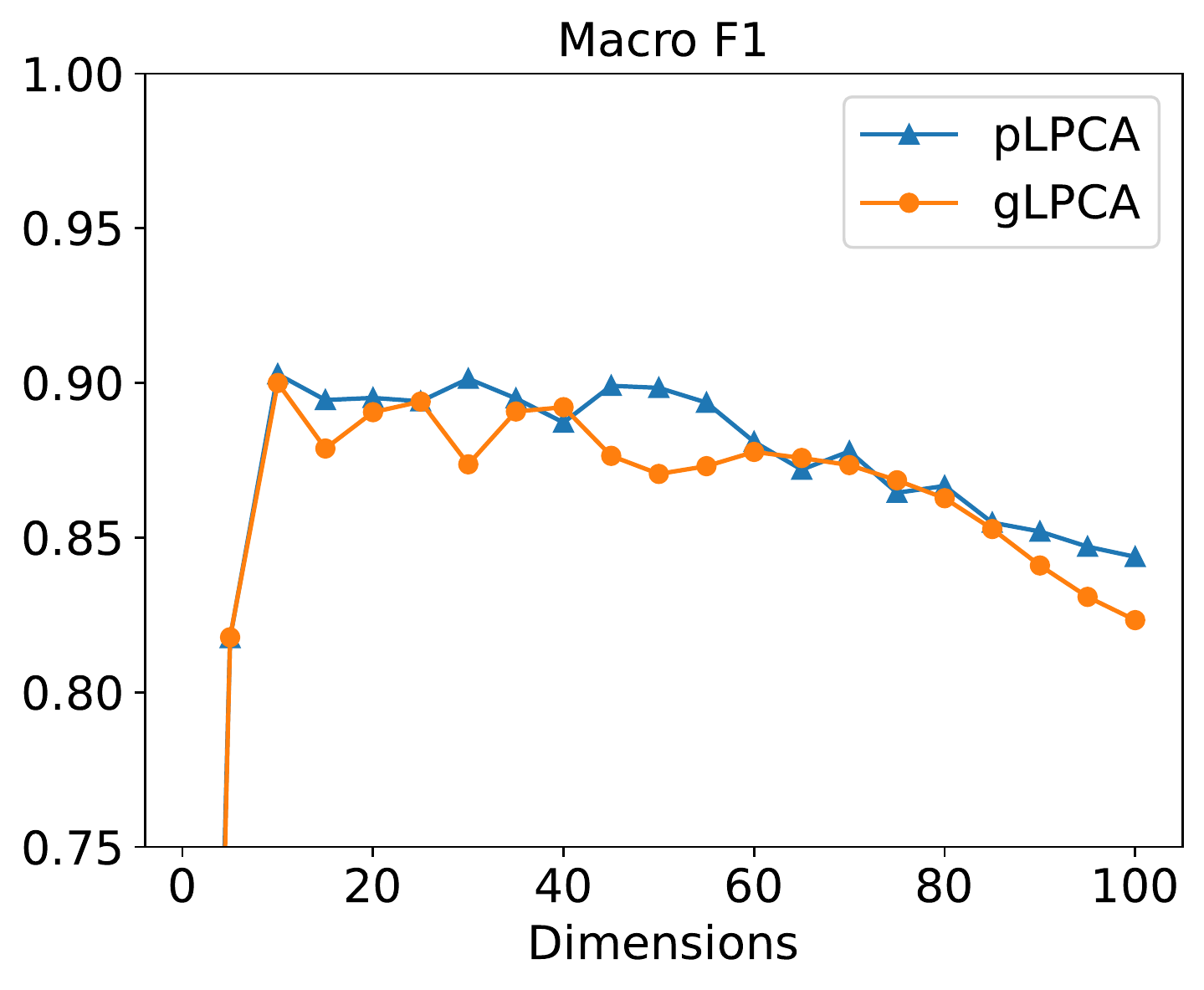}
  \caption{Macro F1 score on the MultiSource dataset across different reduced dimensions (gLPCA vs. pLPCA).}\label{fig:f1fig_pLPCA}
\endminipage
\end{figure}

To provide a visual representation of this performance enhancement, Figure \ref{fig:ACC-pLPCA} illustrates the improvement in Mean Accuracy across each of the tested subspace dimensions ranging from 1 to 100. This visualization clearly demonstrates the superior accuracy achieved by pLPCA compared to gLPCA across almost every reduced dimension. This reinforces the intuitiveness of pLPCA in achieving better accuracy across a wide range of dimensional reductions.

In a similar manner, Figure \ref{fig:f1fig_pLPCA} presents a graphical analysis of the Macro-F1 score. It is evident that, similar to the accuracy results, the Macro-F1 score demonstrates improvement across almost all tested subspace dimensions. This emphasizes the positive impact of incorporating discriminative information, sparseness, and persistent spectral graphs on enhancing performance.

Moving forward, let us explore how the integration of discriminative information, sparseness, and persistent spectral graphs can further enhance performance.

\subsection{Comparison of PLPCA and RLSDSPCA}

After observing the superior performance of pLPCA compared to gLPCA, we conducted a similar study to compare PLPCA and RLSDSPCA, which has been identified as the top-performing method among all PCA-related approaches in the literature \cite{zhang2022enhancing}. Below, we provide a summary of our results on the COAD dataset for both RLSDSPCA and PLPCA.

\begin{table}[H]
    \centering
    \caption{Comparison of PLPCA and RLSDSPCA performance on the COAD dataset.}
    \label{tab: comparison RL vs RPL}
    {
    \begin{tabular}{ ccccccc } 
    \toprule
    Method  & Mean ACC & Mean Macro-REC & Mean Macro-PRE & Mean Macro-F1 & Macro-AUC \\
    \midrule
    \multirow{1}{8em}{RLSDSPCA\cite{zhang2022enhancing}} & 0.9875 & 0.9614 & 0.9504  & 0.9530 &  0.9614\\
    \multirow{1}{8em}{PLPCA} & \textbf{0.9886} & \textbf{0.9680} & \textbf{0.9517} & \textbf{0.9578} & \textbf{0.9680}\\
    \bottomrule
    \end{tabular} 
    }
\end{table}

The results clearly demonstrate that PLPCA outperforms RLSDSPCA in every major category. Specifically, the mean accuracy across all subspaces was 0.9875 for RLSDSPCA and 0.9886 for PLPCA. Additionally, the Macro-F1-score, which highlights the impact of False Negatives, improved from 0.9530 to 0.9578, representing a 0.48\% improvement.

The performance improvement is even more significant when considering the results of the MultiSource dataset, as presented in Table \ref{tab: comparison RL vs RPL}.

\begin{table}[H]
    \centering
    \caption{Comparison of PLPCA and RLSDSPCA performance on the MultiSource dataset.}
    \label{tab: comparison RL vs RPL}
    {
    \begin{tabular}{ ccccccc } 
    \toprule
    Method & Mean ACC & Mean Macro-REC & Mean Macro-PRE & Mean Macro-F1 & Macro-AUC \\
    \midrule
    \multirow{1}{8em}{RLSDSPCA\cite{zhang2022enhancing}}  & 0.9273 & 0.8343 & 0.8972 & 0.8527 &  0.9024\\
    \multirow{1}{8em}{PLPCA} & \textbf{0.9371} & \textbf{0.8393} & \textbf{0.9089} & \textbf{0.8619} & \textbf{0.9065}\\
    \bottomrule
    \end{tabular} 
    }
\end{table}

\begin{figure}[H]
\minipage{0.5\textwidth}
  \includegraphics[width=\linewidth]{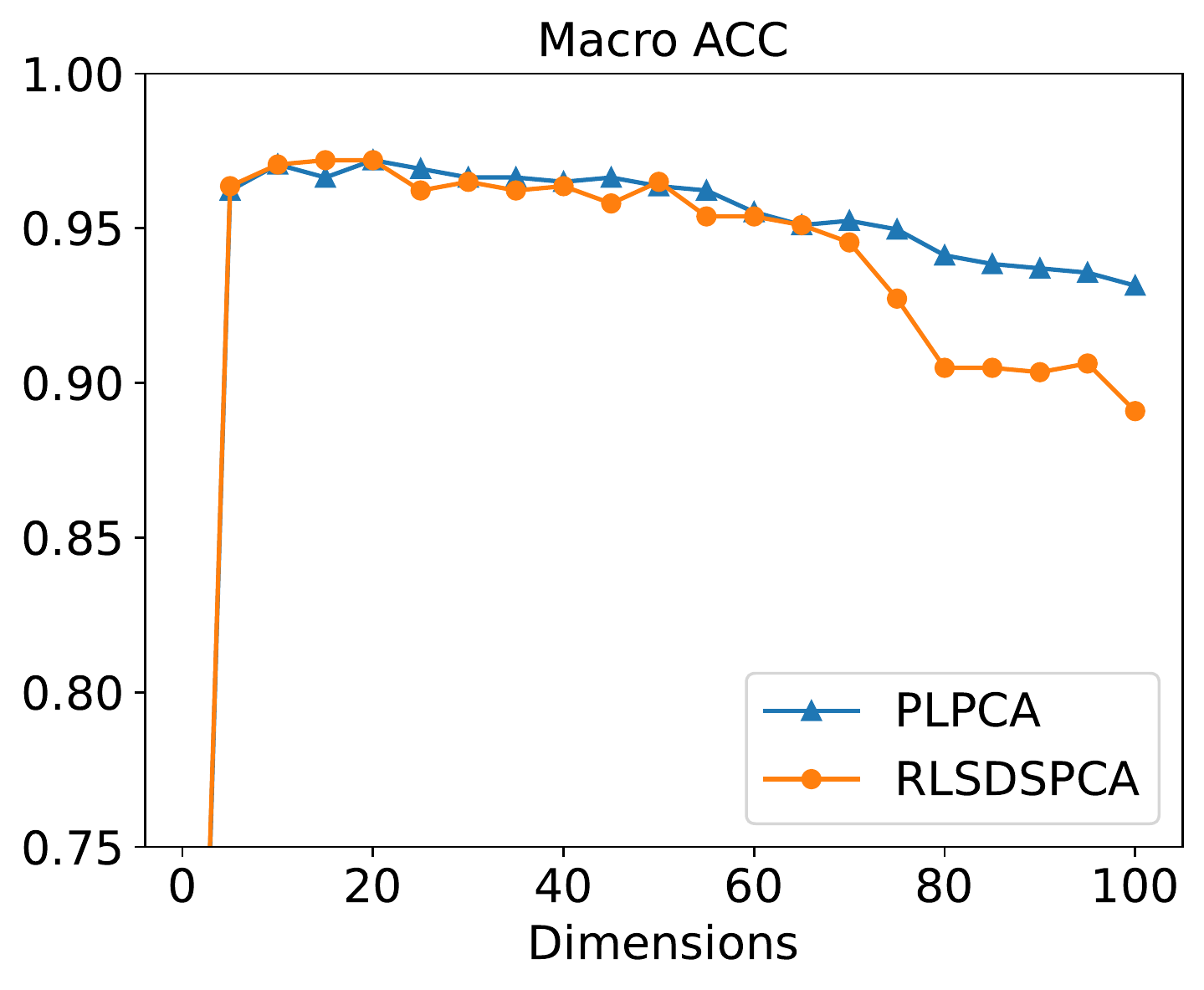}
  \caption{Macro ACC score on the MultiSource dataset across k subspace dimensions (RLSDSPCA vs. PLPCA).}\label{fig:accfig}
\endminipage
\hspace{.5cm}
\minipage{0.5\textwidth}
  \includegraphics[width=\linewidth]{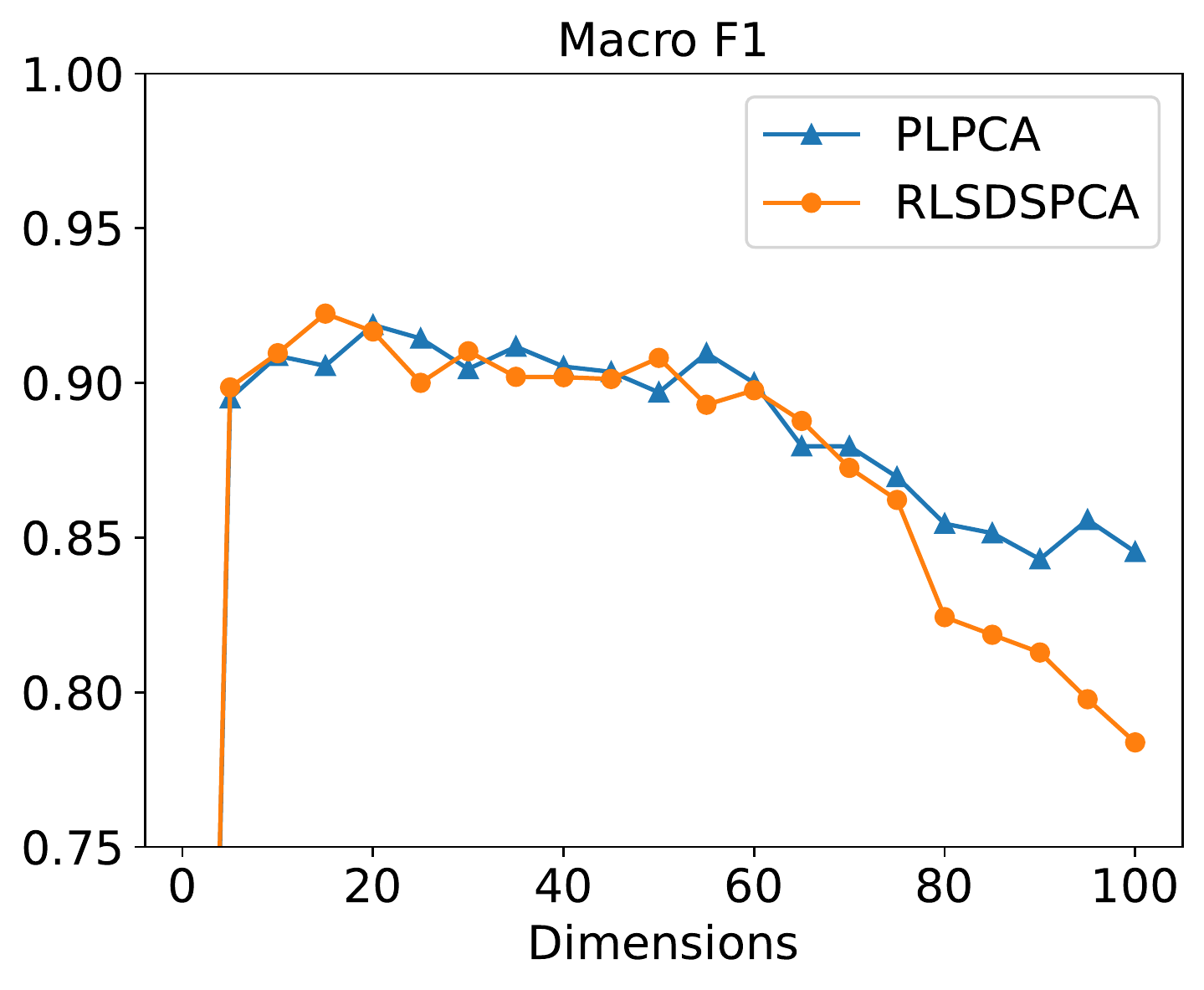}
  \caption{Macro F1 score on the MultiSource dataset across k subspace dimensions (RLSDSPCA vs. PLPCA).}\label{fig:f1fig}
\endminipage
\end{figure}

PLPCA exhibits an improvement in mean accuracy on this benchmark dataset, increasing from 0.9273 to 0.9371, which corresponds to a 0.98\% improvement. Similarly, the F1 score shows an improvement from 0.8527 to 0.8619. To visually demonstrate the comparison between the two methods, Figure \ref{fig:accfig} presents the distribution of performance on the MultiSource dataset across different dimensions for both procedures. This intuitive visualization provides a clear illustration of how the two methods compare in terms of performance.

From the depicted figure, we can observe that while the performance of RLSDSPCA tends to significantly decline as the subspace dimension (k value) increases, this decline is considerably mitigated in the new procedure (PLPCA). As a result, PLPCA exhibits better overall performance. Particularly noteworthy is the even greater improvement in the F1 score, as illustrated in Figure \ref{fig:f1fig}. This demonstrates the ability of PLPCA to consistently outperform the next best method across almost all tested subspace dimensions. This is especially encouraging considering the F1 score's crucial role in tumor classification.

\begin{figure}[H]
\minipage{0.5\textwidth}
  \includegraphics[width=\linewidth]{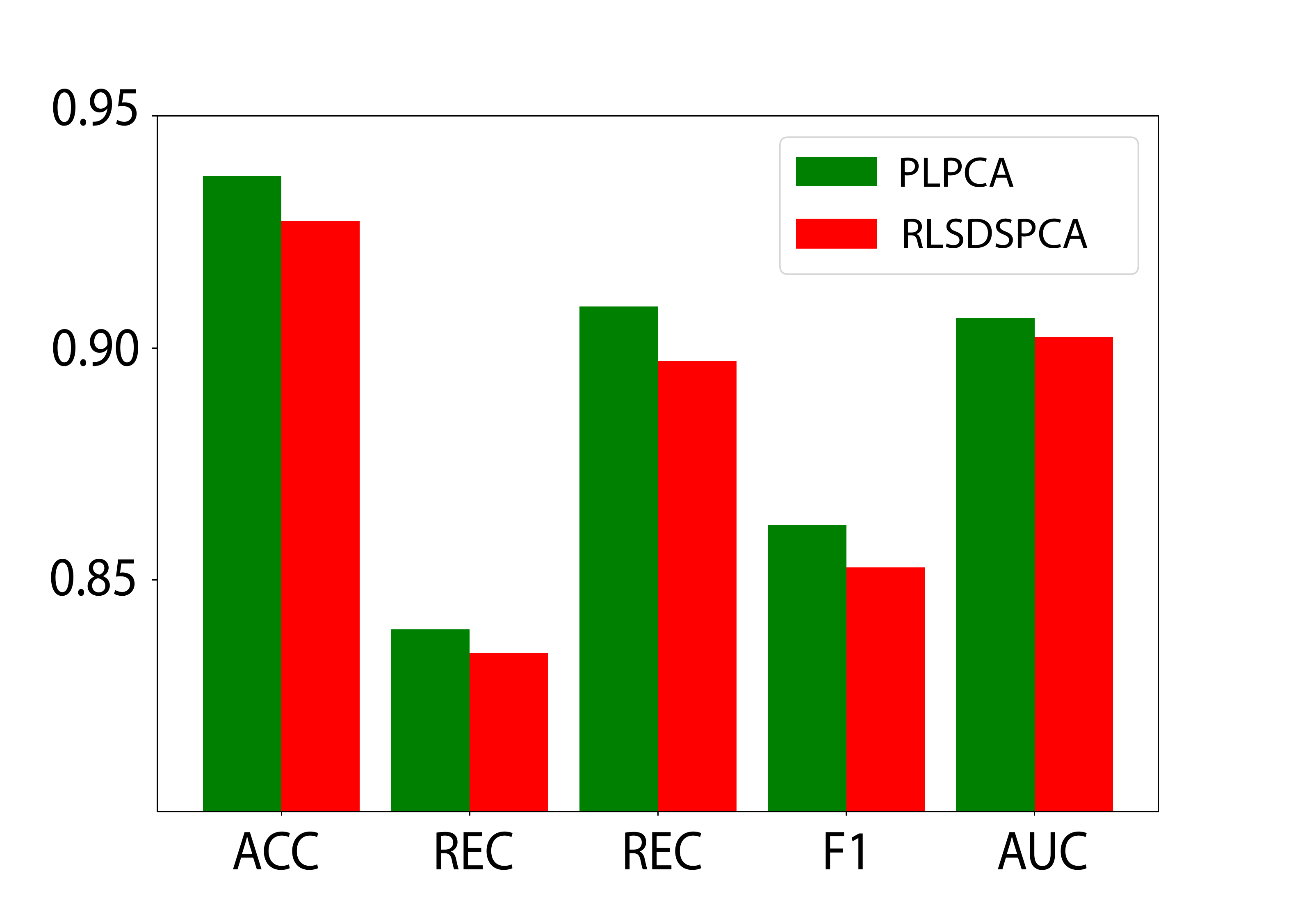}
  \caption{Comparison of each of the evaluation metrics between RLSDSPCA and PLPCA demonstrates classification performance improvement (MultiSource dataset)  }\label{fig:comp1}
\endminipage
\hspace{.5cm}
\minipage{0.5\textwidth}
  \includegraphics[width=\linewidth]{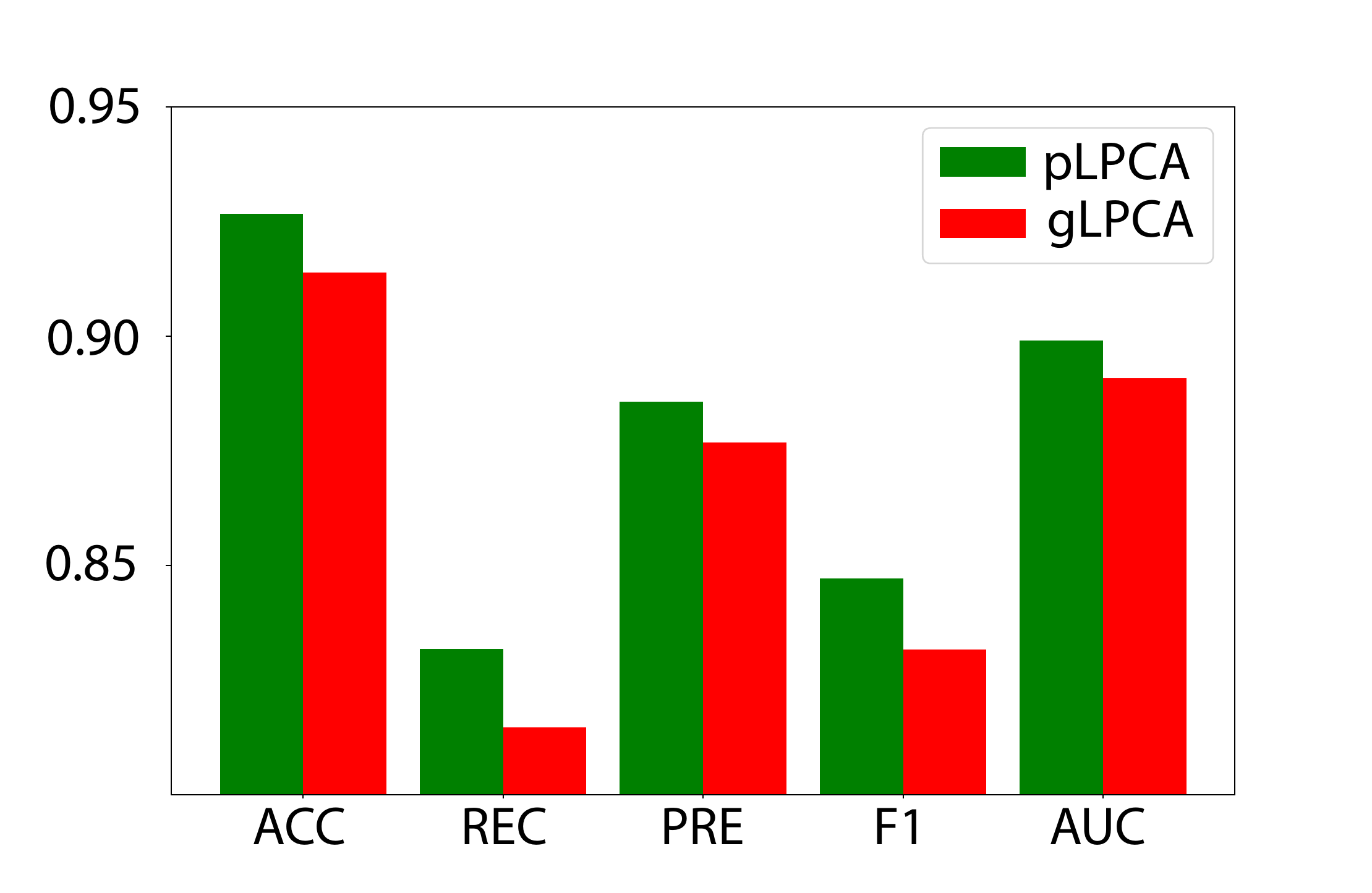}
  \caption{Comparison of each of the evaluation metrics between pLPCA and gLPCA demonstrates classification performance improvement (MultiSource dataset)}\label{fig:comp2}
\endminipage
\end{figure}

To facilitate a more straightforward comparison of the two procedures, we provide a barplot in Figure \ref{fig:comp1}. This barplot allows for a comprehensive evaluation of the relative performances of both procedures across all four evaluation metrics for the MultiSource dataset: Accuracy, Recall, Precision, and F1. We include a similar plot comparing the performances of gLPCA and pLPCA as well in Figure \ref{fig:comp2}.

The figure clearly illustrates the superiority of PLPCA over RLSDSPCA in all the tested metrics, surpassing the previously identified next best method. Notably, the most substantial improvement is observed in the F1 and Recall scores, which are considered particularly important. These findings provide strong evidence of the effectiveness of PLPCA.

To further validate our findings, we  now examine and compare some other PCA methods.

\subsection{Comparisons With Other Methods}

Zhang et al. \cite{zhang2022enhancing} demonstrated that RLSDSPCA achieves superior results compared to other PCA-based approaches. However, after observing how PLPCA enhances classification performance in comparison to RLSDSPCA, it is necessary to further evaluate the performance of our method against other existing approaches.

To validate the performance of PLPCA, we can refer to Table \ref{tab: comparison RPL vs Rest} for the MultiSource dataset, where a comprehensive comparison of different methods is presented.

\begin{table}[H]
    \centering
    \caption{Comparison of PLPCA and other notable methods performance on the MultiSource data}
    \label{tab: comparison RPL vs Rest}
    {
    \begin{tabular}{ lccccc } 
    \toprule
    Method & Mean ACC & Mean Macro-REC & Mean Macro-PRE & Mean Macro-F1 & Macro-AUC \\
    \midrule 
    \multirow{1}{8em}{PLPCA}   & \textbf{0.9371} & \textbf{0.8393} & \textbf{0.9089} & \textbf{0.8619} & \textbf{0.9065}\\
    \multirow{1}{8em}{pLPCA}     & \textbf{0.9267} & \textbf{0.8318} & \textbf{0.8857} & \textbf{0.8471} & \textbf{0.8991}  \\
	   \multirow{1}{8em}{RLSDSPCA\cite{zhang2022enhancing}} &   0.9273    & 0.8343     &  0.8972    & 0.8527 &   0.9024 \\  
	  \multirow{1}{8em}{SDSPCA   \cite{feng2019supervised,zhang2022enhancing}} &  0.9124 & 0.8144 & 0.8917 & 0.8333 & 0.8891\\
	  \multirow{1}{8em}{RgLPCA \cite{jiang2013graph,zhang2022enhancing}} &  0.9197  &  0.8210 &  0.8748  & 0.8353 &  0.8945 \\  
    \multirow{1}{8em}{gLSPCA \cite{feng2019pca,zhang2022enhancing}} & 0.9195    &  0.8148  & 0.8769  &  0.8318 &  0.8910\\ 
    \multirow{1}{8em}{gLPCA \cite{jiang2013graph,zhang2022enhancing}} & 0.9193 & 0.8147 & 0.8768 & 0.8316 & 0.8909 \\
    \multirow{1}{8em}{PCA  \cite{jolliffe2005principal,zhang2022enhancing}} & 0.9108 & 0.7957 & 0.8389 & 0.8025 &   0.8726\\
    \bottomrule 
    \end{tabular} 
    }
    
\end{table} 

This table provides a comprehensive overview of the performance differences between our procedure and other PCA enhancements, including pLPCA. The results clearly demonstrate the significant impact of PLPCA's ability to capture geometrical structure information through persistent spectral graphs, while also incorporating label information and sparseness.

Notably, there is a considerable improvement in Mean Accuracy when transitioning from PCA to PLPCA, with the metric increasing from 0.9108 to 0.9371, representing a 2.63\% improvement. Similarly, the F1 score shows an even greater improvement, increasing from 0.8025 to 0.8619, which corresponds to a remarkable improvement of 5.94\%. These findings underscore the importance of not only capturing geometrical information but also addressing class ambiguities and enforcing sparseness.

Additionally, we can compare the performance of PLPCA with other notable PCA enhancements on the COAD dataset, as depicted in Table \ref{tab: comparison RPL vs Rest}.

\begin{table}[H]
    \centering
    \caption{Comparison of PLPCA and other notable methods performance on the COAD data}
    \label{tab: comparison RPL vs Rest}
    {
    \begin{tabular}{ cccccc } 
    \toprule
    Method & Mean ACC & Mean Macro-REC & Mean Macro-PRE & Mean Macro-F1 & Macro-AUC \\
    \midrule
    \multirow{1}{8em}{PLPCA} & \textbf{0.9886} & \textbf{0.9680} & \textbf{0.9517}  & \textbf{0.9578} & \textbf{0.9680} \\
    \multirow{1}{8em}{pLPCA} &  \textbf{0.9788} & \textbf{0.9450} & \textbf{0.8996} & \textbf{0.9115} & \textbf{0.9450}  \\ 
    \multirow{1}{8em}{RLSDSPCA* \cite{zhang2022enhancing}} & 0.9797& 0.9443& 0.9081& 0.9149& 0.9444\\
    \multirow{1}{8em}{SDSPCA* \cite{feng2019supervised}} & 0.9643 &0.9533  &0.8740  & 0.8918 &0.9533  \\ 
    \multirow{1}{8em}{RgLPCA* \cite{jiang2013graph}} &0.9734    &0.9015   &0.8990   &0.8750  &0.8990   \\  
    \multirow{1}{8em}{gLSPCA* \cite{feng2019pca}} &0.9761    & 0.9250   & 0.8969  & 0.8958 & 0.9250 \\  
    \multirow{1}{8em}{gLPCA* \cite{jiang2013graph}} &0.9756   &0.9429  &0.8841  &0.9002  &0.9429 \\ 
    \multirow{1}{8em}{PCA* \cite{jolliffe2005principal} } &0.9593   &0.8988  &0.8599  &0.8799  &0.8980 \\
    \bottomrule  
    \end{tabular} 
    *Reproduced in the present work
    }
\end{table}

Once again, it is important to highlight the consistently superior performance across all five evaluation metrics, with particular emphasis on the Macro F1 score. The results clearly demonstrate that the PLPCA procedure outperforms other PCA methods by a significant margin.

\begin{figure}[H]
	\centering
	\includegraphics[width = 0.8\textwidth]{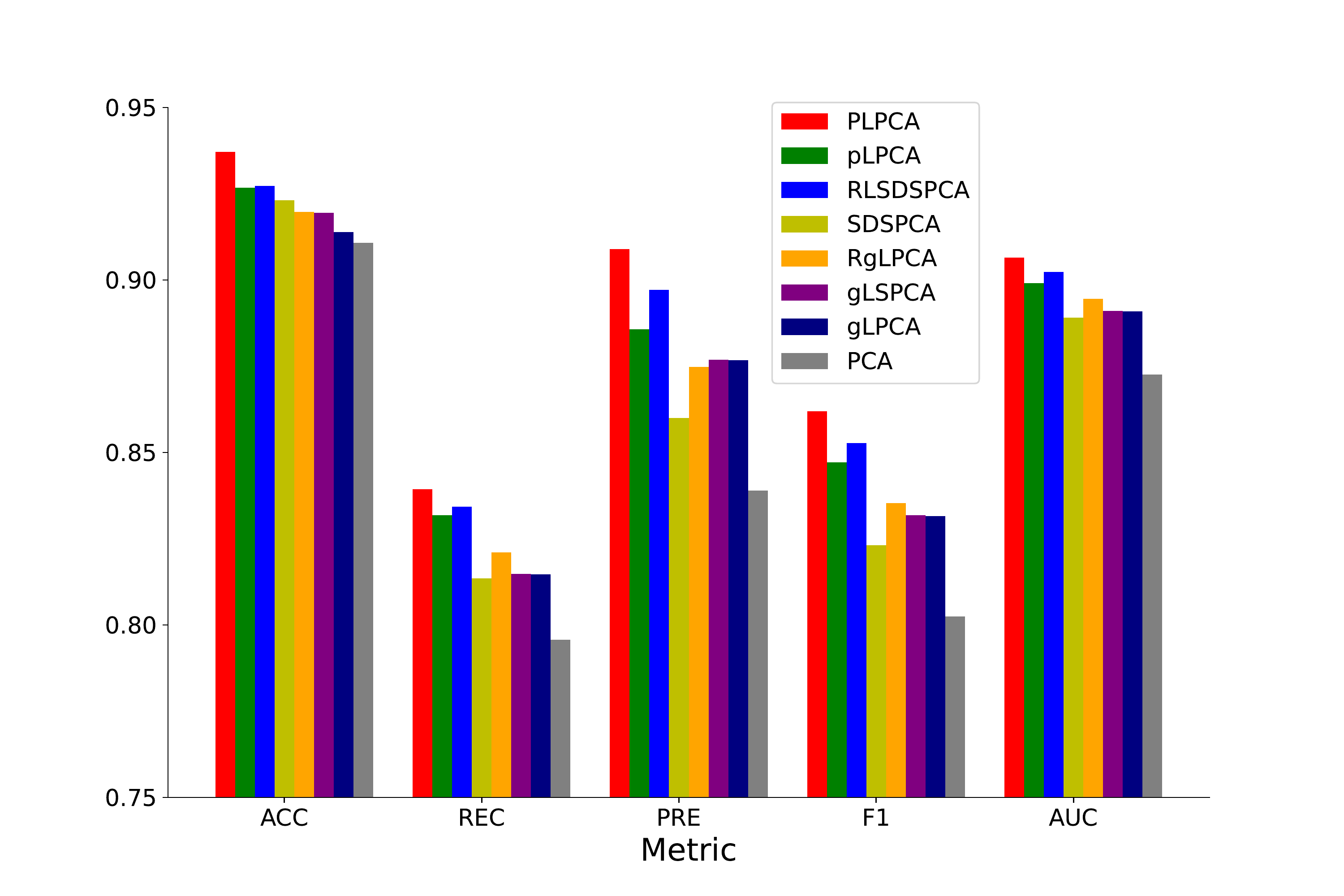}
	\caption{Comparison across five performance metrics for PCA-based methods for the MultiSource dataset.}
	\label{fig:metrics_others_fig}
\end{figure}

To further underscore this point, let us compare our method to traditional PCA. The comparison reveals noteworthy improvements in Mean Accuracy, increasing from 0.9593 to 0.9886, and Mean F1 score, improving from 0.8799 to 0.9578. It is crucial to acknowledge that PLPCA exhibits superior performance across all major evaluation metrics, not just accuracy and F1 score. To visually represent this, we provide a barplot in Figure \ref{fig:metrics_others_fig}, comparing the performance metrics of the mentioned procedures tested on the MultiSource dataset.

From the depicted image, it is evident that PLPCA surpasses other PCA enhancements, including RLSDSPCA, in all aspects, particularly in terms of F1 and Recall.

After confirming the efficacy of our procedure on real gene expression data, we can proceed to evaluate our method on various simulated outlier datasets sourced from the Zhang et al. RLSDSPCA GitHub repository. The objective is to assess whether the enhanced robustness of RLSDSPCA is compromised by the inclusion of persistent spectral graphs.

\subsection{Robustness to Outliers}

\begin{figure}[H]
	\centering
	\includegraphics[width = 0.7\textwidth]{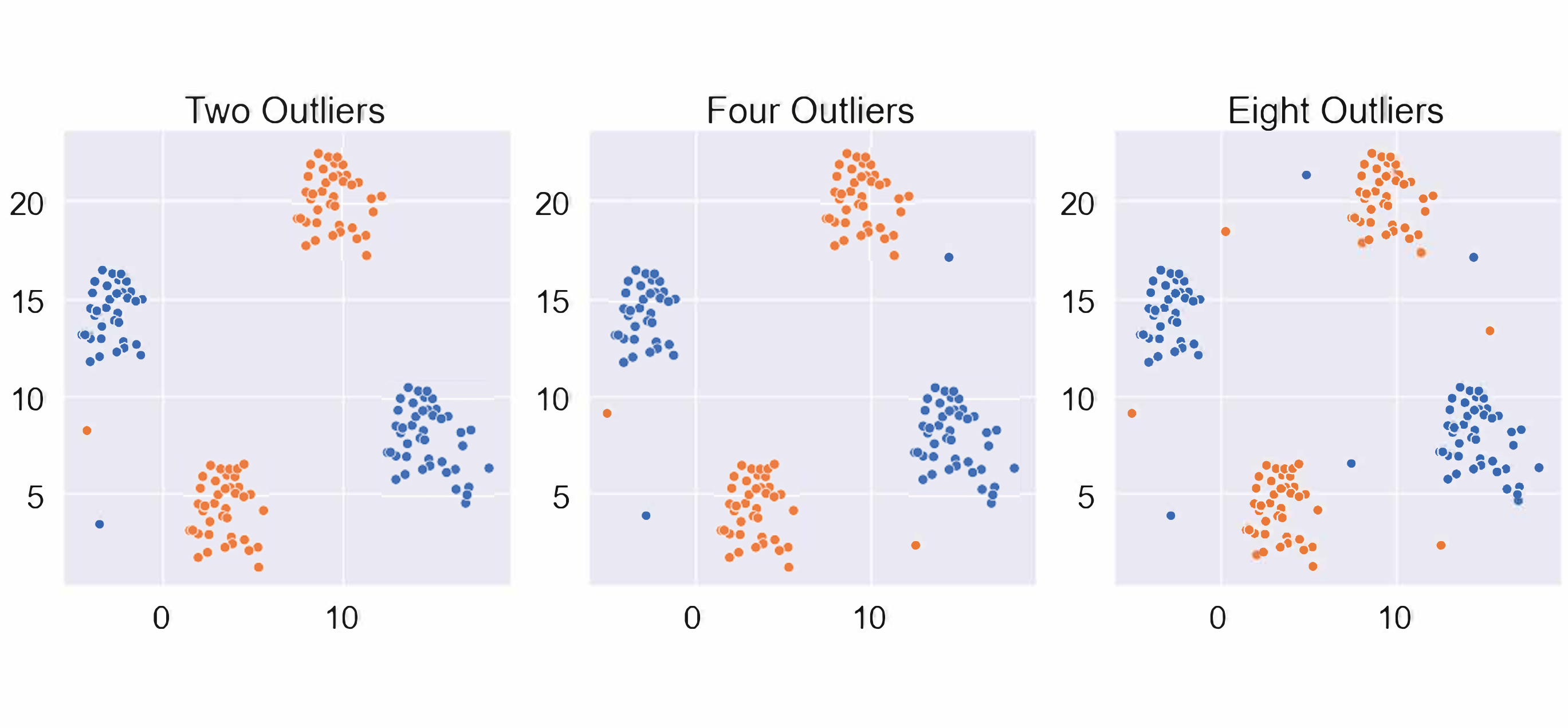}
	\caption{PCA plot of simulated datasets with two, four, and eight outliers for robustness testing. Each dataset has been mapped to k = 2 dimensions for visualization.}
	\label{fig:scatter}
\end{figure}

Earlier studies show that inclusion of $\text{L}_{2,1}$ norm regularization to make the error function robust to outliers \cite{zhang2022enhancing}. Here, we verify the continued effectiveness of this method on the PLPCA procedure by testing on several simulated outlier datasets. The datasets have two, four, and eight outliers respectively. In each case, there are  two classes. We include the PCA plots of each simulated dataset in Figure \ref{fig:scatter}. 

We now summarize the classification results of each algorithm on each dataset in Table \ref{tab: comparison outlier}. 
{\small

\begin{table}[H]
    \centering
    \setlength\tabcolsep{4pt}
    \captionsetup{margin=0.1cm}
    \caption{Verifying robustness on simulated outlier datasets}
    \label{tab: comparison outlier}
    {
    \begin{tabular}{ c|lccccc } 
    \toprule
    Dataset & Method & Mean ACC & Mean Macro-REC & Mean Macro-PRE & Mean Macro-F1  & Macro-AUC\\
    \midrule
    \multirow{6}*{\shortstack{Two Outliers}}  
		    &PLPCA & 1.0 & 1.0 & 1.0 & 1.0 & 1.0\\
        &RLSDSPCA \cite{zhang2022enhancing} & 1.0 & 1.0 & 1.0 & 1.0 &  1.0\\
        &SDSPCA \cite{feng2019supervised,zhang2022enhancing} & 0.9939 & 0.9923 & 0.9952 & 0.9935 &  0.9923   \\
			  &RgLPCA \cite{jiang2013graph,zhang2022enhancing} &0.9939&0.9923&0.9952&0.9935&0.9923\\
        &gLSPCA \cite{feng2019pca,zhang2022enhancing} &0.9939&0.9923& 0.9952&0.9935& 0.9923\\
				&gLPCA \cite{jiang2013graph,zhang2022enhancing} & 0.9939 & 0.9923 & 0.9952 & 0.9935  &  0.9923\\
        &PCA \cite{jolliffe2005principal,zhang2022enhancing}  & 0.9939 & 0.9923 & 0.9952 & 0.9935 & 0.9923 \\
    \midrule
    \multirow{6}*{\shortstack{Four Outliers}}
        & PLPCA & 0.9939 & 0.9952 & 0.9923 & 0.9936 &   0.9952\\
        &RLSDSPCA \cite{zhang2022enhancing} & 0.9939 & 0.9952 & 0.9923 & 0.9936 &   0.9952 \\
        &SDSPCA \cite{feng2019supervised,zhang2022enhancing} & 0.9878 & 0.9889 & 0.9867 & 0.9874 & 0.9889   \\
        &RgLPCA \cite{jiang2013graph,zhang2022enhancing} &0.9818&0.9813& 0.9806&0.9808 & 0.9813\\
        &gLSPCA \cite{feng2019pca,zhang2022enhancing} &0.9878&0.9869&0.9869&0.9869 & 0.9869\\
				&gLPCA \cite{jiang2013graph,zhang2022enhancing} & 0.9818 & 0.9813 & 0.9806 & 0.9808  & 0.9813\\
        &PCA \cite{jolliffe2005principal,zhang2022enhancing} & 0.9818 & 0.9813 & 0.9806 & 0.9808 & 0.9813 \\
    \midrule
    \multirow{6}*{\shortstack{Eight Outliers}}
        &PLPCA & 0.9939 & 0.9947 & 0.9933 & 0.9938 &  0.9947\\
        &RLSDSPCA \cite{zhang2022enhancing} & 0.9939 & 0.9947 & 0.9933 & 0.9938 &  0.9947 \\
        &SDSPCA \cite{feng2019supervised,zhang2022enhancing} & 0.9818 & 0.9823 & 0.9811 & 0.9815 & 0.9823   \\
        &RgLPCA \cite{jiang2013graph,zhang2022enhancing} &0.9757&0.9756&0.9758&0.9754 & 0.9756\\
        &gLSPCA  \cite{feng2019pca,zhang2022enhancing}&0.9818& 0.9823& 0.9811& 0.9815 & 0.9823\\
			  &gLPCA \cite{jiang2013graph,zhang2022enhancing} & 0.9757 & 0.9756 & 0.9758 & 0.9754 &  0.9756 \\
        &PCA \cite{jolliffe2005principal,zhang2022enhancing} & 0.9757 & 0.9756 & 0.9758 & 0.9754 &  0.9756   \\
    \bottomrule
    \end{tabular} 
    }
\end{table}
}
The inclusion of persistent graph regularization in RLSDSPCA does not significantly affect its robustness, indicating that our new method remains robust to outliers to a certain extent. However, it is important to note that an increased number of outliers can still have a negative impact on performance, although this is less problematic for both PLPCA and RLSDSPCA.

\subsection{Residue-Similarity Analysis}

To more effectively visualize our gene expression data after dimensionality reduction, we can generate Residue-Similarity plots for each of the tested datasets \cite{hozumi2022ccp}. We can then compare the results for the gLPCA and pLPCA  models. 

\begin{figure}[H]
	\centering
	\includegraphics[width = .4\textwidth]{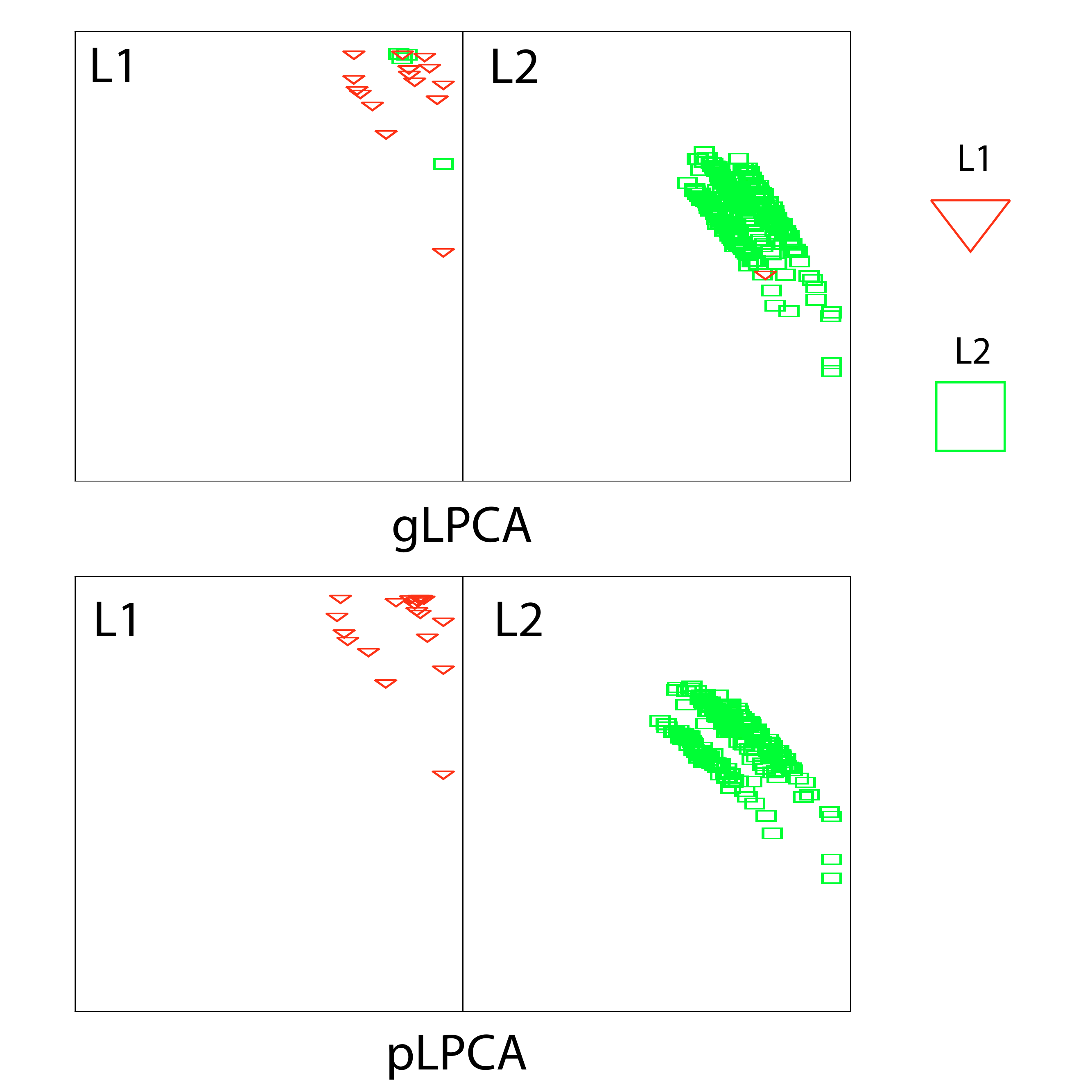}
	\caption{R-S plots of clusters generated from gLPCA and pLPCA-based dimensionality reduction. The $ {x}$-axis is the residual score, and the $ {y}$-axis is the similarity score. Each section corresponds to one cluster and the data were colored according to the predicted labels from KNN on the COAD dataset at $k=100$.}
	\label{fig:rs_COAD}
\end{figure}

\begin{figure}[H]
	\centering
	\includegraphics[width = .8\textwidth]{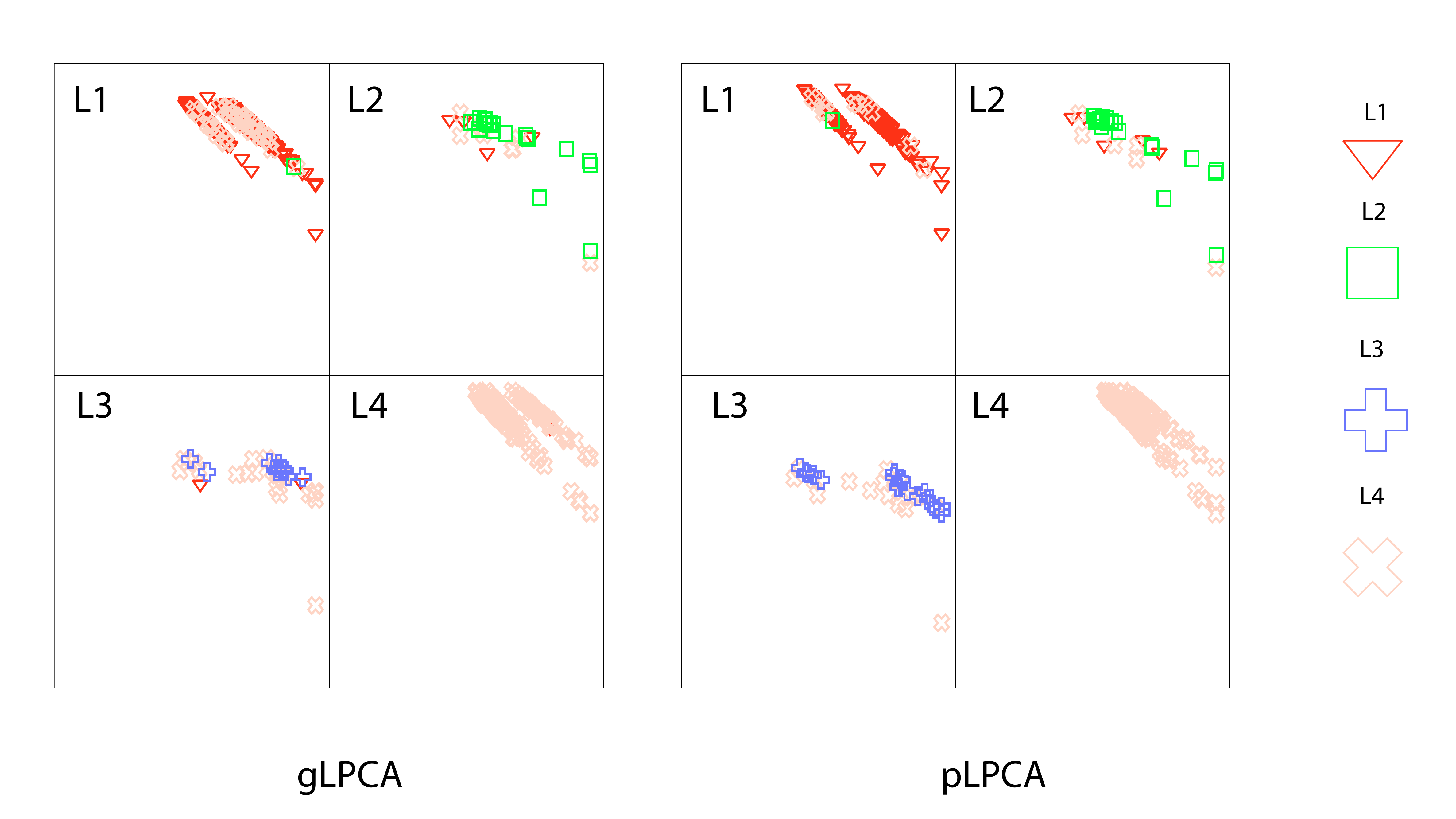}
	\caption{R-S plots of clusters generated from gLPCA and pLPCA-based dimensionality reduction. The $ {x}$-axis is the residual score, and the $ {y}$-axis is the similarity score. Each section corresponds to one cluster and the data were colored according to the predicted labels from KNN on the MultiSource dataset at $k=65$.} 
	\label{fig:rs_MS}
\end{figure}

Figures \ref{fig:rs_COAD} and \ref{fig:rs_MS} compare the classification accuracy on the COAD and MultiSource datasets between dimensionality reductions using gLPCA and pLPCA. The chosen subspace dimensions for visualization were $k = 100$ and $k = 65$ respectively. These results, along with Figures \ref{fig:ACC-pLPCA} and \ref{fig:f1fig_pLPCA} showcase the ability of persistent Laplacian-regularized PCA to outperform graph Laplcian-based PCA. In particular, we note the poor performance of gLPCA-based dimensionality reduction for classifying Labels 1 and 3 in the MultiSource dataset, and the improvement seen when using pLPCA instead. 

\section{Conclusion}

As DNA sequencing technologies have advanced in timeliness and cost, they have greatly expanded our understanding of pathogenic genes responsible for the development and progression of different cancers. These insights have led to the identification of diagnostic biomarkers and therapeutic targets. Given the need for dimensionality reduction to effectively analyze this data, it is crucial to maximize the accurate representation of the data. In this regard, we propose a novel method called Persistent Laplacian-enhanced PCA (PLPCA). This method incorporates robustness, label information, and sparsity, while also improving the capture of geometrical structure information using techniques derived from persistent topological Laplacian theory \cite{wang2020persistent}.

Our extensive computational results demonstrate that by incorporating persistent topological regularization in the RLSDSPCA procedure, we achieve the highest level of classification performance after dimensionality reduction compared to previous methods. While the inclusion of a graph Laplacian contributes to capturing geometrical structure information, its analysis is limited to a single-scale Laplacian. To overcome this limitation, we generate a sequence of topological Laplacians through filtration, providing a more comprehensive multiscale perspective of the data and enabling us to emphasize features at important scales. Alternatively, we also achieve the similar superb  results by incorporating the PL regularization  to the original PCA approach. This method, called pLPCA does not depend on the availability of data labels.

Despite progress made by our proposed method, there is still ample room for improvement. 
First, it is interesting to examibne the role of   higher-dimensional Laplacians in dimensionality reduction.  Additionally, further analysis is needed to evaluate the performance of this procedure for feature selection compared to other methods. Previous studies, including Zhang et al. \cite{zhang2022enhancing}, have described a feature selection procedure that assumes linear relationships among genes, which may not be optimal \cite{huerta2008analysis}. It would be advantageous to explore more sophisticated feature selection techniques that account for the nonlinear relationships among genes. Moreover, integrating our new dimensionality reduction procedure into these methods could lead to further performance improvements \cite{kiselev2017sc3, ren2019sscc}. Additionally, understanding the role and significance of the selected genes in driving or correlating with different cancer incidences is an important area for future research. Both aspects require continued efforts in the fields of mathematics and biology.

\section{Data and Model Availability}
The model and data used in this analysis is publicly available at the \href{https://github.com/seanfcottrell/PLPCA}{PLPCA} GitHub repository: \\
https://github.com/seanfcottrell/PLPCA

\section*{Acknowledgements}
This work was supported in part by NIH grant  GM126189, NSF grants DMS-2052983,  DMS-1761320, and IIS-1900473,  NASA grant 80NSSC21M0023,  Michigan Economic Development Corporation, MSU Foundation,  Bristol-Myers Squibb 65109, and Pfizer.

\section{Conflicts of Interests}
The authors declare no competing interests. 

\clearpage 

\bibliographystyle{unsrt}
\bibliography{refs}

\end{document}